\documentclass[sn-mathphys]{sn-jnl}

\usepackage{amsmath,paralist,amssymb,amsfonts,epsfig,float,color,comment,subfigure,caption,textcomp,mathtools}
\everymath{\displaystyle}

\jyear{2021}%

\theoremstyle{thmstyleone}%
\newtheorem{theorem}{Theorem}
\theoremstyle{thmstyletwo}%

\theoremstyle{thmstylethree}%
\theoremstyle{thmstylethree}%
\newtheorem{lem}{Lemma}%

\begin{document}

\title[Mathematical modeling of multilane vehicular traffic by discrete kinetic theory]{Mathematical modeling of multilane vehicular traffic by discrete kinetic theory approach}


\author{\fnm{Mohamed} \sur{Zagour}}\email{m.zagour@insa.ueuromed.org}

\affil{\orgdiv{Euromed Research Center, INSA Euro-Mediterranean,} \orgname{Euromed University of Fes, Morocco}}


\abstract{This paper deals with the modeling and numerical simulations of multilane vehicular traffic according to the discrete kinetic theory approach. The nonlinear additive interactions and external actions such as tollgates as well traffic signs are considered and modeled. The well-posedness of the related Cauchy problem for the spatially homogeneous case has been proved by using Banach fixed-point theory. Numerical simulations are carried out to validate the improved model in the cases of spatially homogeneous and inhomogeneous problems.}




\maketitle

\section{Introduction}
As it is known, the mathematical approach  to vehicular traffic  modeling is developed on three representation scales: microscopic, macroscopic, and kinetic. However, the critical analysis proposed in \cite{[BDS11]} states that none of those scale approaches is totally satisfactory. Precisely,  each one presents technical and conceptual advantages and disadvantages. Consequently, a multiscale approach is a necessity to obtain a detailed description of the dynamics of vehicles on the road. The author in \cite{[Do2]} compliments this point of view, where a hybrid model is proposed with detailed modeling of the dynamics of the micro-systems which is implemented into a macroscopic hyperbolic equation. An important reference to clarify research activity in the field is the survey of \cite{[HEL01]} on the physics and modeling of multi-particle systems. While, the critical paper by \cite{[DAG95]} provides some drawbacks of the driver-vehicle micro-system, where interactions can even modify the behavior of the driver and his/her ability is conditioned by the local flow conditions. Next, this paper has generated various discussions and reactions to account for the aforementioned criticisms~\cite{[BDF12],[HJ09],[PSTV15],[PSTV17]}. We mention that the approach of \cite{[BDF12]} has been further developed by various papers, for instance by a multiscale continuous velocity and activity variables for vehicular traffic flow \cite{[CNZ19]}.

This paper deals with the mathematical modeling of vehicular traffic flow along a multilane road on the basis of the kinetic theory. Precisely, we refer to Reference \cite{[BG09]}, where the authors proposed a discrete kinetic model with the main features:
\begin{enumerate}
	\item The approach is developed at the kinetic scale where the modeling of the interacting vehicles is considered at the microscopic scale with a binary interactions;
	\item The velocity variable is assumed to be discrete to overtake the drawback of vehicles number might not large enough to assure continuity of the probability distributions over the micro-states; 
	\item The quality of the road conditions is taken into account by an additional parameter which takes values in the interval $[0, 1]$, where the extremes of the interval correspond to worst and best conditions respectively.
\end{enumerate}

The proposing discrete model can be considered as a deep revisiting of the proposed  model in \cite{[BG09]}. Indeed, the aforesaid paper consider only the numerical simulations of the spatially homogeneous problem where the asymptotic property in time has been shown. In this paper, we consider both the spatially homogeneous and inhomogeneous problems. In the first case, the well-posedness of the related Cauchy problem has been proved according to Banach fixed-point theory, and the Kerner's fundamental diagrams has been reproduced. In the second problem, we show the numerical results of the emerging clusters phenomena. In the modeling part, we deal with the following topics: $i)$ nonlinear additive interactions (rather than binary interactions) between vehicles accounting on $ii)$ perceived  density (rather than real one), $iii)$  function depending on the space variable modeling the road conditions (rather than a parameter), and $iv)$ dynamics under external actions such as the presence of tollgates. The modeling is inspired by the books of \cite{[KER04],[PH71]} which reports a detailed interpretation of the physics traffic.

Let us comment now on some references in the literature about the modeling of multilane traffic flow. One can find the macroscopic approach in \cite{[MBY84]} while \cite{[MSR10]} reports an interesting survey of different models for lane changing. Also, a multilane model is analyzed in \cite{[GKR03]} focusing on the total vehicle density across all lanes, while a modeling traffic vehicular with conservation laws is proposed in \cite{[GR20]}. On the basis of kinetic theory, microscopic models are studied in \cite{[KW99a],[KW99b]} and a macroscopic model derived from individual behavior is given in \cite{[HG97]}. Recently, the authors in \cite{[HMV18]} extended the analysis to the second-order Aw-Rascle model and a hybrid stochastic kinetic model, respectively. Note that the approach in \cite{[CC06]} is more similar to the analysis presented in \cite{[HN19]} but with a different source term.

The rest of this paper is organized as follows: Section \ref{Sect2} deals with a concise description of multilane traffic flow and derives a new mathematical structure suitable to include the aforementioned features in addition to those already included in \cite{[BDF12]}. Section \ref{Sect3} is devoted to the derivation of a specific model by inserting a possible form of the terms that appeared in the proposed mathematical structure, that is to say, the perceived density, the encounter rate, table of games, the weight function and the external force. Section \ref{Sect4} presents a qualitative analysis of the spatially homogeneous problem. Moreover, Kerner's fundamental diagrams \cite{[KER04]} are reproduced and the trend to equilibrium (asymptotic in time) is shown. Finally, in Section \ref{Sect5} the numerical simulations of the spatially inhomogeneous problem are provided to capture emerging behaviors.
\section{Mathematical representation and structure}\label{Sect2}
This section is organized through two subsections: the first one deals with the description of traffic flow along a multilane and the second one is devoted to the design of mathematical frameworks for nonlinearly additive interactions along with multilane traffic flow.
\subsection{Representation of traffic flow along a multilane}
This subsection deals with a discrete kinetic description of traffic flow along a multilane road and introduces the classical macroscopic quantities.\\
Herein, all the dependent variables describing the traffic flow are dimensionless and normalized to values in the interval $[0,1]$.
\begin{itemize}
	\item[\textbullet] $X$ is the length of the road;
	\item[\textbullet] $V_M $ is the maximum velocity attained by an isolated fast vehicle when moving in free flow conditions;
	\item[\textbullet] $\tilde{\rho}_M$ is the maximum density of vehicles corresponding to bumper to bumper vehicle distance.
\end{itemize}
Using the above quantities, the following dimensionless variables are identified:
\begin{itemize}
	\item[\textbullet] $x$ is the position referred to length of road $X$;
	\item[\textbullet] $t$ is the time normalized by means of $\frac{X}{V_M}$;
	\item[\textbullet] $v=\frac{V}{V_M}$ is the velocity referred to $V_M$;
	\item[\textbullet] $\rho=\tilde{\rho}/\tilde{\rho}_M$ is the density of vehicles referred to the maximum one.
\end{itemize}
In this paper, we adopt the discrete kinetic theory method by introducing the  uniform grids of the velocity variable taking the following forms
\begin{equation}
	0=v_1<v_2<\dots<v_i<\dots<v_n=1, \qquad i=1,\dots,n.
\end{equation}

Under these hypotheses, the evolution in time and space of the multilane traffic flow is described by introducing in each lane a statistical distribution function as a linear combination of Dirac distributions in variables $v$ and $u$ as follows: 
$$f^\ell(t,x,v)=\sum_{i=1}^nf_{i}^\ell(t,x)\otimes\delta_i(v-v_i):\mathbb{R}^+\times[0,1]\times[0,1]\rightarrow\mathbb{R}^+,$$
where  $f_{i}^\ell(t,x)=f^\ell(t,x,v_i)$ for $\ell=1,\cdots,L$ with $L$ is the total number of lanes. Thus, macroscopic quantities can derived in each lane as follow
\begin{itemize}
	\item[\textbullet] The density in the $\ell$-lane is
	$$\rho_\ell(t,x)=\sum_{i=1}^nf_{i}^\ell(t,x); $$
	\item[\textbullet] The flux in the $\ell$-lane is
	$$q^\ell(t,x)=\sum_{i=1}^nv_if_{i}^\ell(t,x);$$
	\item[\textbullet] The average velocity in the $\ell$-lane is
	$$U^\ell(t,x)=\frac{q^\ell(t,x)}{\rho^\ell(t,x)}=\frac{1}{\rho_\ell(t,x)}\sum_{i=1}^nv_if_{i}^\ell(t,x);$$
	\item[\textbullet] The variance of the velocity in the $\ell$-lane is 
	$$\Theta^\ell(t,x)=\frac{1}{\rho_\ell(t,x)}\sum_{i=1}^n\big(v_i-U^\ell(t,x)\big)^2f_{i}^\ell(t,x).$$
\end{itemize}
The global density and the global flux are obtained by summing the contributions of all lanes
$$\rho(t,x)=\sum_{\ell=1}^L\rho_\ell(t,x),\qquad q(t,x)=\sum_{\ell=1}^Lq^\ell(t,x).$$
Moreover, the global average velocity and variance are given by
$$U(t,x)= \frac{1}{L}\sum_{\ell=1}^LU^\ell(t,x)\qquad \Theta(t,x)=\sum_{\ell=1}^L\Theta^\ell(t,x).$$
\subsection{A mathematical structure toward modeling}
This subsection reports about the derivation of a specific structure in the case of nonlinearly additive interactions involving the micro-systems, generating the evolution of the distribution function $f^\ell_{i}$. In general, the kinetic approach is such that interactions are modeled by table of games. Thus, three types of vehicles are involved:

\noindent \textbullet $\;$ \emph{Test particle}: It represents the whole system. The related distribution function is  $f^\ell_{i}=f^\ell(t,x,v_i)$;\\
\noindent \textbullet $\;$ \emph{Candidate particle}: It changes its current state to that of the test particle as a consequence of an interaction. The related distribution function is
$f^r_{h}=f^r(t,x,v_h)$;\\
\noindent \textbullet $\;$ \emph{Field particle}: It interacts with test and candidate vehicles. The related distribution function is  $f^r_{p}=f^r(t,x^*,v_p)$.\\
Our structure is obtained by a balance of vehicles in the elementary volume of the space of the microscopic state which includes position  and velocity. This balance of vehicles includes the free transport term, the dynamics of nonlinear additive interactions, and the trend to the velocity imposed by the external actions. We mention that the dynamics of nonlinear additive interactions include a gain term of vehicles that enter in the
aforementioned elementary volume and a loss term of vehicles that leave it. Thus, the resulting structure can be written, at a formal level, as follows:
\begin{equation}
	\displaystyle\partial_t f_{i}^\ell(t,x)+v_i\,\partial_x f_{i}^\ell(t,x)
	=  \mathcal{J}_{i}^\ell[\textbf{f},\textbf{f}](t,x)+\mathcal{T}_{i}^\ell[\textbf{f}](t,x),
\end{equation}
where $f^\ell_{i} = f^\ell(t, x, v_i),\;\textbf{f}=(f_1,\cdots,n)$ and $v_i\,\partial_xf^\ell_{i}
$ is the free flow transport term, while $\mathcal{J}_{i}^\ell$ and $\mathcal{T}^\ell_{i}$ correspond, respectively, to nonlinear additive interactions and interaction with external actions. Thus, the evolution of the distribution function in each lane is given by
{\small\begin{eqnarray}\label{Struct}
	\begin{array}{l}
		\displaystyle\partial_t f_{i}^\ell(t,x)+v_i\,\partial_x f_{i}^\ell(t,x)
		=  \mathcal{G}_{i}^\ell[\textbf{f},\textbf{f}](t,x)-f_{i}^\ell(t,x)\,\mathcal{L}^\ell[\textbf{f}](t,x)+\mathcal{T}_{i}^\ell[\textbf{f}](t,x)\\ 
		{}\\ 
		\displaystyle\hskip 0.1cm  =\sum_{r=1}^L\sum_{h,p=1}^n\int_{\Omega_v}\eta^r[\rho^\star_r(t,x^*)]\mathcal{A}_{hp,r}^{i,\ell}[\rho_1^\star,\dots,\rho_L^\star;\alpha](t,x^*)f_{h}^{r}(t,x)f_{p}^{r}(t,x^*)\omega^r(x,x^*)dx^* \\
		{}\\ 
		\displaystyle\hskip.5cm-f_{i}^\ell(t,x)\int_{\Omega_v}\eta^\ell[\rho^\star_\ell(t,x^*)]\rho_\ell(t,x^*)\,\omega^\ell(x,x^*)dx^*+\mu^\ell[\rho^\star_\ell]\Big(f_{e}^\ell(x;v_e(x))-f_{i}^\ell(t,x)\Big).
	\end{array}
\end{eqnarray}}
In model \eqref{Struct}, $\Omega_v=[x,x+\xi_v]$ represents the visibility zone where $\xi_v$ is the visibility length on front of the vehicle, that depends on the quality of the road-environment, that is to say on $\alpha = \alpha(x)$; $\eta^\ell[\rho_\ell(t,x^*),x]$ is the encounter rate, which depends on the probability distributions by means of the perceived density $\rho^\star_\ell$ in the $\ell-$lane; $\mathcal{A}_{h,p,r}^{i,j,\ell}$ defines the table of games, which denotes the probability density that the candidate particle falls into the state of the test particle after an interaction with a field particle; $\omega^\ell(x,x^*)$ represents the weight function in each lane $\ell$; and $\mu^\ell[\rho^\star_\ell]$  models the intensity of the action, which increases with $\rho_\ell^\star$, while $v_{e}(x)$ is the speed imposed by the external action in the $\ell-$lane.\\
Herein, we shall consider admissible interactions that generate transitions only in continuous lanes. In details, the sets of admissible lane transitions are as follows
{\small\begin{equation}
	I_r=
	\left\{
	\begin{array}{ll}
		I_r=\{1,2\},& r=1;\\
		I_r=\{r-1,r,r+1\}, & r\neq 1,L; \\
		I_r=\{L-1,L\}, &r=L.
	\end{array}
	\right.
\end{equation} }

\section{From mathematical structure to model}\label{Sect3}
This section deals with the derivation of a specific model of multilane vehicular traffic by inserting into the aforesaid mathematical structure (\ref{Struct}) models of interactions at the microscopic scale. This goal is pursued by looking at the modeling of the interaction terms that characterize such structure, that is to say, $\eta^\ell, \mathcal{A}_{hp,r}^{i,\ell}, \omega^\ell,\mu^\ell$ and $f_e^\ell$, such that a good agreement with empirical data, concerning both the Kerner's fundamental diagram and the emerging behaviors in unsteady flow conditions, is given.
\subsection{Modeling of interaction domain and perceived density}
The vehicles have a visibility zone denoted by $\Omega_v = \Omega_v(x) = [x, x + \xi]$, where $\xi$ is the visibility length on front of the particle, that depends on the quality of the road-environment conditions modeled by $\alpha = \alpha(x)$.  In more detail, we assume that $\xi = \alpha \, X$, where $X << \xi$ is the visibility length in the case of best quality of the road ($\alpha = 1$). Note that one can consider a sensitive zone $\Omega_s=[x,x+\xi_s]$, necessary to perceive the flow conditions in $\Omega_s$. We refer the interested reader to the paper \cite{[CNZ19]} for more details.\\
The concept of perceived density was introduced in~\cite{[DEA99]}, where it was suggested that this quantity is greater (smaller) than the real one whenever positive (negative) density gradients appear. The following expression is considered
{\small\begin{equation}
	\rho^\star_\ell[f]  =  \rho_\ell + \frac{\partial_x \rho_\ell}{\sqrt{1 + (\partial_x \rho_\ell)^2}}\,\Big((\frac{1}{L}- \rho_\ell)\, H(\partial_x \rho_\ell) + \rho_\ell \, H(- \partial_x \rho_\ell)\Big),
\end{equation}}
where  $H(\cdot)$ is the Heaviside  function $H(\cdot \geq 0) = 1$, while $H(\cdot < 0) = 0$. Thus,  the perceived density, positive gradients increase the value of $\rho^\star_\ell$ from $\rho_\ell$ to the maximum admissible value $\rho_\ell = \frac{1}{L}$, while negative gradients decrease it from $\rho_\ell$ to the lowest admissible value $\rho_\ell = 0$ such that
$$
\partial_x \rho_\ell \to +\infty \Rightarrow \rho^\star_\ell \to \frac{1}{L}, \quad \partial_x  \rho_\ell = 0 \Rightarrow \rho^\star_\ell = \rho_\ell,  \quad \partial_x  \rho_\ell \to - \infty \Rightarrow \rho^\star_\ell \to 0.
$$
\subsection{Modeling of the encounter rate and the weight function}
The encounter rate $\eta^\ell[\rho^\star_\ell(t,x)]$ gives the number of interactions per unit time among the micro-systems in each lane. We assume that this term grows with the local perceived density starting from a minimum value corresponding to driving in vacuum conditions $\eta^0$. The following expression is considered
$$\eta^\ell[\rho^\star_\ell(t,x)]=\eta^0\big(1+\gamma_\eta^\ell L\rho^\star_\ell(t,x)\big),$$
where $\gamma_\eta^\ell$ is the growth coefficient and $\rho^\star_\ell$ is the perceived density in each lane.\\
The weight function $\omega^\ell(x,x^*)$ is assumed to be the same in all lanes and it satisfies the following requirement:
$$\omega^\ell(x,x^*)\geq 0, \qquad \int_{\Omega_v}\omega^\ell(x,x^*)dx^*=1, \qquad \forall x^*\in {\Omega_v}.$$

\subsection{Modeling of the table of games}
We present and discuss a possible form for the table of games $A_{h,p,r}^{i,\ell}$ which gives the probability that a micro-system in the lane $r$ with velocity $v_h$ reaches the velocity $v_i$ in the lane $\ell$, after an interaction with a micro-system traveling at velocity $v_p$ in the lane $r$. It satisfies the following 
\begin{equation}\label{tablB}
	\mathcal{A}_{hp,r}^{i,\ell}\geq 0,\quad \sum_{\ell=1}^L\sum_{i=1}^n\mathcal{A}_{hp,r}^{i,\ell}=1,
	\quad \forall h,p=1,\cdots,n,\quad \forall r=1,\cdots,L
\end{equation} 

Herein, we improve the table of games proposed by \cite{[BG09]} by taking into account the nonlinear additive interactions rather than binary ones, perceived density rather than the real one, and space function modeling the road conditions rather than a parameter.\\

We consider that the candidate particle with velocity $v_h$ and assume that it interacts with a field particle with velocity $v_p$ in the same lane. Thus, the cases $v_h > v_p$, $v_h < v_p$, $v_h = v_p$, are analysed separately. Herein, we omit the dependence of the road conditions and the perceived density functions on the space variable $x$. \\
\textbf{ I. Interaction with a faster particle $(v_h < v_p)$.}
In this case, the candidate particle is encountering a faster field particle in the $r$-lane. Then the candidate particle either maintains its lane or it changes it to the right continuous one, which means $I_r=\{r-1,r\}$. In the following we distinguish the case $r=\{2,\cdots,L\}$ from the case $r=1$. 
\subparagraph*{a) \underline{Interaction in the case $h\neq1,$ $r\neq1$}}
As a result of the above assumptions, it follows that $I_r=\{r-1,r\}$. Therefore,
{\small\begin{equation}\label{3.3}
	\mathcal{A}_{hp,r}^{i,\ell=r}=
	\left\{
	\begin{array}{ll}
		0, & i=1,\cdots,h-1, \\
		(1-\alpha)\big(1-L\rho^\star_r(1-L\rho^\star_{r-1})\big), & i=h,\\
		\displaystyle\alpha \frac{1}{(i-h)}\frac{1}{\sum_{z=h+1}^{p}(z-h)^{-1}}\big(1-L\rho^\star_r(1-L\rho^\star_{r-1})\big),&i=h+1,\cdots,p.
	\end{array}
	\right.
\end{equation}
\begin{equation}\label{3.4}
	\mathcal{A}_{hp,r}^{i,\ell=r-1}=
	\left\{
	\begin{array}{ll}
		\displaystyle (1-\alpha )(h-i)\frac{1}{\sum_{z=1}^{h-1}(h-z)}L\rho^\star_r(1-L\rho^\star_{r-1}), & i=1,\cdots,h-1, \\
		\displaystyle\alpha  L\rho^\star_r(1-L\rho^\star_{r-1}),&i=h,\\
		0,&i=h+1,\cdots,p.
	\end{array}
	\right.
\end{equation}}
Note that Eq. (\ref{3.3}) gives the probability that the candidate particle maintains its lane. It can maintain its current speeds or accelerate. In the latter case, it can reach a new velocity $v_i \in\{v_{h+1},\cdots,v_p\}$ with a probability that depends not only on the road conditions and on the perceived density but also on the distance between the velocity classes involved. Eq. (\ref{3.4}) yields the probability that the candidate particle changes its lane to the right lane.
\subparagraph*{b) \underline{Interaction in the case $h=1$ and $r\neq1$}}
In this case, $I_r=\{r-1,r\}$ and the field particle stops or the candidate particle accelerates in the same lane. Therefore,
{\small\begin{equation}
	\mathcal{A}_{1p,r}^{i,\ell=r}=
	\left\{
	\begin{array}{ll}
		L\rho^\star_r, & i=1, \\
		\displaystyle\frac{1}{(i-1)}\frac{1}{\sum_{z=2}^{p}(z-1)^{-1}}(1-L\rho^\star_r), & i=2,\cdots,p.
	\end{array}
	\right.
\end{equation}}
\subparagraph*{c) \underline{Interaction in the case $h=1,\cdots,p-1$ and  $r=1$}}
In this case, the candidate particle is in the slowest lane $(I_r=1)$. We assume that when it interacts with a faster field particle, it can only maintain its current lane with the same velocity or accelerating. In the latter case, it can reach a new velocity $v_i \in\{v_{h+1},\cdots,v_p\}$ with a probability that depends not only on the road conditions and on the perceived density but also on the distance between the velocity classes involved. Consequently,
{\small\begin{equation}
	\mathcal{A}_{hp,1}^{i,\ell=1}=
	\left\{
	\begin{array}{ll}
		0,&i=1,\cdots,h-1, \\
		\displaystyle1- \alpha (1-L\rho^\star_1), & i=h, \\
		\displaystyle\alpha \frac{1}{(i-h)}\frac{1}{\sum_{i=h+1}^{p-1}(i-h)^{-1}}\big(1-L\rho^\star_1\big) ,&i=h+1,\cdots,p.
	\end{array}
	\right.
\end{equation}}
\textbf{II. Interaction with a faster particle $(v_h > v_p)$.}
In this case, the candidate particle is encountering a slower field particle in the $r-$lane. Thus, the candidate particle either maintains its lane or it changes to the left continuous one, which means $I_r = \{r, r+1\}$. In what follows, we distinguish the case $r = \{1,\cdots,L-1\}$ where the candidate particle can change its lane, from the case $r = L$ in which it can only maintain the current lane.
\subparagraph*{a) \underline{Interaction in the case $h\neq n$ and $r\neq L$}}

In this case $I_c=\{r,r+1\}$. Therefore,

{\small\begin{equation}\label{3.7}
	\mathcal{A}_{hp,r}^{i,\ell=r}=
	\left\{
	\begin{array}{ll}
		L\rho^\star_{r+1}, & i=p, \\
		0,&otherwise.
	\end{array}
	\right.
\end{equation}
\begin{equation}\label{3.8}
	\mathcal{A}_{hp,r}^{i,\ell=r+1}=
	\left\{
	\begin{array}{ll}
		\displaystyle(1-\alpha )\frac{h-i}{\sum_{i=p+1}^{h}{(h-z)}}(1-L\rho^\star_{r+1}), & i=p+1,\cdots,h, \\
		\displaystyle\alpha \frac{1}{(i-h)}\frac{1}{\sum_{z=h+1}^{n-1}(z-h)^{-1}}(1-L\rho^\star_{r+1}),&i=h+1,\cdots,n-1.
	\end{array}
	\right.
\end{equation}}
Eq. (\ref{3.7}) yields the probability that the candidate particle maintains its lane when interacting with a slower field particle traveling in the same  $r-$lane. The probability depends only on the perceived density in the $(r+1)-$lane. Eq. (\ref{3.8}) yields the probability that the candidate particle changes its position to the left lane. Here, the probabilities depend on $\alpha$, $\rho^\star_{r+1}$, the difference between the velocity classes involved.
\subparagraph*{b) \underline{Interaction in the case $h= n$ and $r\neq L$}}
As a consequence of the above assumptions, it follows that $I_r=\{r,r+1\}$. In this case, the velocity ${v_n}$ of the candidate particle is the maximum allowed. It maintains the current lane reducing its velocity to that of the field particle or it changes to the left lane. While in the left lane, the candidate particle maintains its velocity or brakes depending on the local traffic density. Thus,
{\small\begin{equation}
	\mathcal{A}_{np,r}^{i,\ell=r}=
	\left\{
	\begin{array}{ll}
		L\rho^\star_{r+1}, & i=p, \\
		0,&otherwise.
	\end{array}
	\right.
\end{equation}
\begin{equation}
	\mathcal{A}_{np,r}^{i,\ell=r+1}=
	\left\{
	\begin{array}{ll}
		\displaystyle\frac{n-i}{\sum_{z=p+1}^n(n-z)}(1-L\rho^\star_{r+1}), & i=p+1,\cdots,n, \\
		0,& otherwise.
	\end{array}
	\right.
\end{equation}}
\subparagraph*{c) \underline{Interaction in the case $h=p+1,\cdots,n$ and $r= L$}}
The candidate particle is in the fastest lane. When it interacts with a slower field particle, we assume that it is obligated to travel with the velocity $v_p$. Thus,
{\small\begin{equation}
	\mathcal{A}_{hp,L}^{i,\ell=L}=
	\left\{
	\begin{array}{ll}
		1, & i=p, \\
		0,&otherwise.
	\end{array}
	\right.
\end{equation}}
\textbf{III. Interaction with an equally faster particle $(v_h=v_p)$.}
In this case, the candidate particle and the field particle travel with the same speed in the same $r-$lane. We assume that the candidate particle maintains its lane or it can change both to the right lane or in the left one, which means $I_r=\{r-1,r,r 1\}$. In what follows, we distinguish the case in which the candidate particle can change both to the right and the left lanes from those where it can change only to the right or only to the left lane. 
\subparagraph*{a) \underline{Interaction in the case $h\neq 1,n$ and $r\neq 1,L$}} As a consequence of the above assumptions, it follows that $I_r=\{r-1,r,r+1\}$. Thus,
{\small\begin{equation}\label{3.12}
	\mathcal{A}_{hp,r}^{i,\ell=r}=
	\left\{
	\begin{array}{ll}
		\displaystyle\frac{1}{2}(1-\alpha )(h-i)\frac{1}{\sum_{z=1}^{h-1}(h-z)}L\rho^\star_r(L\rho^\star_{r-1}+L\rho^\star_{r+1}),& i=1,\cdots,h-1, \\
		\displaystyle\frac{1}{2}(1-\alpha )(1-L\rho^\star_r)(L\rho^\star_{r-1}+L\rho^\star_{r+1}),&i=h,\\
		\displaystyle\frac{1}{2}\alpha \frac{1}{(i-h)}\frac{1}{\sum_{z=h+1}^{n}(z-h)^{-1}}(L\rho^\star_{r-1}+L\rho^\star_{r+1}),&i=h+1,\cdots,n.
	\end{array}
	\right.
\end{equation}
\begin{equation}\label{3.13}
	\mathcal{A}_{hp,r}^{i,\ell=r-1}=
	\left\{
	\begin{array}{ll}
		\displaystyle\frac{1}{2}(1-\alpha )(h-i)\frac{1}{\sum_{z=1}^{h-1}(h-z)}L\rho^\star_r(1-L\rho^\star_{r-1}),& i=1,\cdots,h-1, \\
		\displaystyle\frac{1}{2}\alpha (1-L\rho^\star_{r-1}),& i=h,\\
		0,&i=h+1,\cdots,n.
	\end{array}
	\right.
\end{equation}
\begin{equation}\label{3.14}
	\mathcal{A}_{hp,r}^{i,\ell=r+1}=
	\left\{
	\begin{array}{ll}
		0,& i=1,\cdots,h-1, \\
		\displaystyle\frac{1}{2}(1-\alpha )(1-L\rho^\star_{r+1}),& i=h,\\
		\displaystyle\frac{1}{2}\alpha \frac{1}{(i-h)}\frac{1}{\sum_{z=h+1}^{n}(z-h)^{-1}}(1-L\rho^\star_{r+1}),&i=h+1,\cdots,n.
	\end{array}
	\right.
\end{equation}}
Eq. (\ref{3.12}) yields the probability that the candidate particle maintains its current lane $r$, in this case, the candidate particle can decelerate or accelerate or maintains its current speed. All the above probabilities depends on the road conditions $\alpha$, the densities $\rho_{r-1}$, $\rho_{r}$, $\rho_{r+1}$, the difference between the velocity classes involved. Eq. (\ref{3.13}) yields the probability that the candidate particle changes its lane to the right lane, in this case the candidate particle cannot accelerate. Eq. (\ref{3.14}) yields the probability that the candidate particle changes its lane to the left lane, in this case, the candidate particle cannot brake.

\subparagraph*{b) \underline{Interaction in the case $h\neq 1,n$ and $r=1$}} 
In this case, the interacting vehicles are in the slowest lane and $I_r=\{1,2\}$. Consequently,
{\small\begin{equation}\label{3.15}
	\mathcal{A}_{hp,1}^{i,\ell=1}=
	\left\{
	\begin{array}{ll}
		\displaystyle(1-\alpha )(h-i)\frac{1}{\sum_{z=1}^{h-1}(h-z)}L\rho^\star_1L\rho^\star_2,& i=1,\cdots,h-1,\\
		\displaystyle(1-\alpha )(1-L\rho^\star_{1})L\rho^\star_2,& i=h,\\
		\displaystyle\alpha \frac{1}{(i-h)}\frac{1}{\sum_{z=h+1}^{n}(z-h)^{-1}}L\rho^\star_{2},&i=h+1,\cdots,n.
	\end{array}
	\right.
\end{equation}
\begin{equation}\label{3.16}
	\mathcal{A}_{hp,1}^{i,\ell=2}=
	\left\{
	\begin{array}{ll}
		0,& i=1,\cdots,h-1,\\
		\displaystyle(1-\alpha )(1-L\rho^\star_{2}),& i=h,\\
		\displaystyle\alpha \frac{1}{(i-h)}\frac{1}{\sum_{z=h+1}^{n-1}(z-h)^{-1}}(1-L\rho^\star_{2}),&i=h+1,\cdots,n.
	\end{array}
	\right.
\end{equation}}
Eq. (\ref{3.15}) yields the probability that the candidate particle maintains its current lane, depending on the two traffic densities, on the parameter $\alpha$, on the difference between the velocity classes involved. Eq. (\ref{3.16}) yields the probability that the candidate particle changes its lane to the left one, in this case it cannot brake.
\subparagraph*{c) \underline{Interaction in the case $h\neq 1,n$ and $r=L$}} 
In this case, the interacting vehicles are in the fastest lane and  $I_r=\{L-1,L\}$. The candidate particle or it changes to the right one, depending on the two traffic densities, on the parameter $\alpha$, on the difference between the velocity classes involved. If it changes its lane, it cannot accelerate. Thus,
{\small\begin{equation}\label{3.17}
	\mathcal{A}_{hp,L}^{i,\ell=L}=
	\left\{
	\begin{array}{ll}
		\displaystyle(1-\alpha )(h-i)\frac{1}{\sum_{z=1}^{h-1}(h-z)}L\rho^\star_{L-1}L\rho^\star_L,& i=1,\cdots,h-1, \\
		\displaystyle(1-\alpha )(1-L\rho^\star_{L})L\rho^\star_{L-1},& i=h,\\
		\displaystyle\alpha \frac{1}{(i-h)}\frac{1}{\sum_{z=h+1}^{n}(z-h)^{-1}}L\rho^\star_{L-1},&i=h+1,\cdots,n.
	\end{array}
	\right.
\end{equation}
\begin{equation}\label{3.18}
	\mathcal{A}_{hp,L}^{i,\ell=L-1}=
	\left\{
	\begin{array}{ll}
		\displaystyle(1-\alpha )(h-i)\frac{1}{\sum_{z=1}^{h-1}(h-z)}(1-L\rho^\star_{L-1}),&i=1,\cdots,h-1, \\
		\displaystyle\alpha (1-L\rho^\star_{L-1}),& i=h,\\
		0,&i=h+1,\cdots,n.
	\end{array}
	\right.
\end{equation}}
Eq. (\ref{3.17}) and (\ref{3.18}) yields the probabilities that the candidate particle maintains its current lane or it changes to the right one, depending on the two traffic densities, on the $\alpha$, on the difference between the velocity classes involved. If it changes its lane it cannot accelerate. 
\subparagraph*{d) \underline{Interaction in the case $h=1$ and $r\neq 1,L$}} 

As a result of the above assumptions, it follows that $I_r=\{r-1,r,r+1\}$. Therefore 
{\small\begin{equation}\label{3.19}
	\mathcal{A}_{11,r}^{i,\ell=r}=
	\left\{
	\begin{array}{ll}
		\displaystyle 1-\alpha L\rho_r\big(1-\frac{1}{2}(L\rho^\star_{r-1}+L\rho^\star_{r+1})\big),& i=1, \\
		0,& i=2,\cdots,n.
	\end{array}
	\right.
\end{equation}
\begin{equation}\label{3.20}
	\mathcal{A}_{11,r}^{i,\ell=r-1}=
	\left\{
	\begin{array}{ll}
		0,& i=1, \\
		\displaystyle\frac{1}{2}\alpha \frac{1}{(i-1)}\frac{1}{\sum_{i=2}^{n}(i-2)^{-1}}L\rho^\star_r(1-L\rho^\star_{r-1}),& i=2,\cdots,n.
	\end{array}
	\right.
\end{equation}
\begin{equation}\label{3.21}
	\mathcal{A}_{11,r}^{i,\ell=r+1}=
	\left\{
	\begin{array}{ll}
		0,& i=1, \\
		\displaystyle\frac{1}{2}\alpha \frac{1}{(i-1)}\frac{1}{\sum_{i=2}^{n}(i-2)^{-1}}L\rho^\star_r(1-L\rho^\star_{r+1}),& i=2,\cdots,n.
	\end{array}
	\right.
\end{equation}}
Eq. (\ref{3.19}) yields the probability that the candidate particle remains in the same lane, and it maintains its velocity. Eqs. (\ref{3.20}) (\ref{3.21}) yields the probabilities that the candidate particle changes its lane to the right lane or the left one, depending on the traffic densities, on the quality of the road, on the difference between the velocity classes involved and on the activity $u_k$. In this case, it can only accelerate and it cannot brake.
\subparagraph*{e) \underline{Interaction in the case $h=1$ and $r=1$}} 

In this case, the two interacting vehicles are in the slowest lane and $I_r =\{1,2\}$. Consequently
{\small\begin{equation}\label{3.22}
	\mathcal{A}_{11,1}^{i,\ell=1}=
	\left\{
	\begin{array}{ll}
		\displaystyle 1-\alpha L\rho^\star_1(1-L\rho^\star_{2}),& i=1, \\
		0,& i=2,\cdots,n.
	\end{array}
	\right.
\end{equation}
\begin{equation}\label{3.23}
	\mathcal{A}_{11,1}^{i,\ell=2}=
	\left\{
	\begin{array}{ll}
		0,& i=1, \\
		\displaystyle\alpha \frac{1}{(i-1)}\frac{1}{\sum_{z=2}^{n}{(z-1)^{-1}}}L\rho^\star_1(1-L\rho^\star_{2}),& i=2,\cdots,n.
	\end{array}
	\right.
\end{equation}}
Eq. (\ref{3.22}) yields the probability that the candidate particle remains in the same lane and only maintains its velocity $v_1$, depending on the traffic densities, on the quality of the road, on the difference between the velocity classes involved. Eq. (\ref{3.23}) yields the probability that the candidate particle changes its lane to right one, in this case it can accelerate.
\subparagraph*{f) \underline{Interaction in the case $h=1$ and $r=L$}} As a consequence of the above assumptions, it follows that $I_r=\{L-1,L\}$. In this case, the two interacting vehicles stop in the fastest lane. Thus,
{\small\begin{equation}\label{3.24}
	\mathcal{A}_{11,L}^{i,\ell=L}=
	\left\{
	\begin{array}{ll}
		\displaystyle 1-\alpha L\rho^\star_L(1-L\rho^\star_{L-1}),& i=1, \\
		0,& i=2,\cdots,n.
	\end{array}
	\right.
\end{equation}
\begin{equation}\label{3.25}
	\mathcal{A}_{11,L}^{i,\ell=L-1}=
	\left\{
	\begin{array}{ll}
		0,& i=1, \\
		\displaystyle\alpha \frac{1}{(i-1)}\frac{1}{\sum_{z=2}^{n}(z-1)^{-1}}L\rho^\star_L(1-L\rho^\star_{L-1}),& i=2,\cdots,n.
	\end{array}
	\right.
\end{equation}}
Eq. (\ref{3.24}) yields the probabilities that the candidate particle remains in the same lane and maintains its current velocity $v_1$, depending on the traffic densities, on the quality of the road,  on the difference between the velocity classes involved. Eq. (\ref{3.25}) yields the probabilities that the candidate particle changes its lane to right one accelerating, depending on the traffic densities, on the quality of the road,  on the difference between the velocity classes involved.
\subparagraph*{g) \underline{Interaction in the case $h=n$ and $r\neq 1,L$}}  In this case, $I_r=\{r-1,r,r+1\}$ and the two interacting vehicles travel at the fastest allowed velocity $v_n=1$ and they are both in a central lane. It is considered that the candidate particle maintains its current lane with the same or a slower velocity or it can change going to a continuous one. Consequently,
{\small\begin{equation}
	\mathcal{A}_{nn,r}^{i,\ell=r}=
	\left\{
	\begin{array}{ll}
		\displaystyle\frac{1}{2}(n-i)\frac{1}{\sum_{z=1}^{n-1}(n-z)}(L\rho^\star_{r-1}+L\rho^\star_{r+1})L\rho_r,& i=1,\cdots,n-1, \\
		\displaystyle\frac{1}{2}(L\rho^\star_{r-1}+L\rho_{r+1})(1-L\rho_r),& i=n.
	\end{array}
	\right.
\end{equation}
\begin{equation}
	\mathcal{A}_{nn,r}^{i,\ell=r-1}=
	\left\{
	\begin{array}{ll}
		\displaystyle\frac{1}{2}(1-\alpha )(n-i)\frac{1}{\sum_{z=1}^{n-1}(n-z)}(1-L\rho^\star_{r-1}),& i=1,\cdots,n-1, \\
		\displaystyle\frac{1}{2}\alpha (1-L\rho^\star_{r-1}),& i=n.
	\end{array}
	\right.
\end{equation}
\begin{equation}
	\mathcal{A}_{nn,r}^{i,\ell=r+1}=
	\left\{
	\begin{array}{ll}
		0,& i=1,\cdots,n-1, \\
		\displaystyle\frac{1}{2}(1-L\rho^\star_{r+1}),& i=n.
	\end{array}
	\right.
\end{equation}}
\subparagraph*{h) \underline{Interaction in the case $h=n$ and $r=1$}}  In this case, $I_r =\{1,2\}$ and the two interacting vehicles are both in the slowest lane with the maximum allowed speed. The candidate particle maintains its current lane traveling with the same or lower velocity independence of the traffic density $\rho_1$ or it changes to the left lane maintaining its velocity. Consequently,
{\small\begin{equation}
	\mathcal{A}_{nn,1}^{i,\ell=1}=
	\left\{
	\begin{array}{ll}
		\displaystyle(n-i)\frac{1}{\sum_{z=1}^{n-1}(n-z)}L\rho^\star_1L\rho^\star_2,& i=1,\cdots,n-1,\\
		\displaystyle(1-L\rho^\star_1)L\rho^\star_2,& i=n.
	\end{array}
	\right.
\end{equation}
\begin{equation}
	\mathcal{A}_{nn,1}^{i,\ell=2}=
	\left\{
	\begin{array}{ll}
		0,& i=1,\cdots,n-1, \\
		1-L\rho^\star_2,& i=n.
	\end{array}
	\right.
\end{equation}}
\subparagraph*{i) \underline{Interaction in the case $h=n$ and $r=L$}} As a result of the above assumptions, it follows the $I_r=\{L-1,L\}$. In this case, the two interacting vehicles are both in the fastest lane with the maximum allowed speed $v_n$. The candidate particle maintains its current lane traveling with the same or lower velocity independence of the traffic density $\rho_L$ or alternatively it changes to the right lane maintaining its velocity or braking. Consequently,
{\small\begin{equation}
	\mathcal{A}_{nn,L}^{i,\ell=L}=
	\left\{
	\begin{array}{ll}
		\displaystyle(n-i)\frac{1}{\sum_{z=1}^{n-1}(n-z)}L\rho_{L-1}L\rho^\star_L,& i=1,\cdots,n-1, \\
		\displaystyle(1-L\rho^\star_L)L\rho^\star_{L-1},& i=n.
	\end{array}
	\right.
\end{equation}
\begin{equation}
	\mathcal{A}_{nn,L}^{i,L=L-1}=
	\left\{
	\begin{array}{ll}
		\displaystyle(1-\alpha )(n-i)\frac{1}{\sum_{z=1}^{n-1}(n-z)}(1-L\rho^\star_{L-1}),& i=1,\cdots,n-1, \\
		\alpha (1-L\rho^\star_{L-1}),& i=n.
	\end{array}
	\right.
\end{equation}}
\subsection{Modeling of external actions}
The last term needed to model is the external action which indicates a prescribed speed as it occurs, for instance, in the presence of tollgates and traffic signs.
The structure of this term is reported in Eq.~(\ref{Struct}), where  $f_e^\ell$ is a given function of the prescribed velocity $v_e= v_e(x)$. Thus, from Eq. (\ref{Struct}) it only requires to model the intense of the action. We propose the following form 
\begin{equation}\label{mu}
	\mu^\ell[\rho^\star_\ell]= \eta^0 \, (1 + \gamma_\mu^\ell \, L\rho_\ell^\star),
\end{equation}
where $\gamma_\mu^\ell$ is the growth coefficient and $\rho^\star$ is the perceived density.
\section{The spatially homogeneous problem}\label{Sect4} In this section, we address the theoretical and computational analysis of the spatially homogeneous problem in which the distribution function $f^\ell$ is independent of the space variable $x$. We define it as follows
$$f^\ell=f^\ell(t,v)=\sum_{i=1}^nf_{i}^\ell(t)\otimes\delta_i(v-v_i):\mathbb{R}^+\times[0,1]\rightarrow\mathbb{R}^+.$$
Consequently, it results that
$$\partial_x f_{i}^\ell=0, \qquad \forall i=1,\cdots,n, \quad \forall \ell=1,\cdots,L.$$
This implies $\rho^\star_\ell=\rho_\ell$ for $ \; \ell=1,\cdots,L$.
Taking into account the above hypotheses, the final form of model \eqref{Struct} in the spatially homogeneous case is reduced to the following Cauchy problem
{\small\begin{equation}\label{SptialP}
	\left\{
	\begin{array}{l}
		\displaystyle\frac{df_{i}^\ell}{dt}=\sum_{r\in I_c}\eta^r[\rho_r(t)]\sum_{h,p=1}^n\mathcal{A}_{hp,r}^{i,\ell}f_{h}^rf_{p}^r-\eta^\ell[\rho_\ell(t)]f_{i}^\ell\rho_\ell(t), \\
		\displaystyle f_{i}^\ell(0)=f_{i,0}^{\ell},
	\end{array}
	\right.
\end{equation}}
where $$\mathcal{A}_{h,p,r}^{i,\ell}=\mathcal{A}_{hk,pq,r}^{i,\ell}[v_h\rightarrow v_i,I_r\rightarrow I_\ell\mid,v_h,v_p,I_r,\rho_1(t),\dots,\rho_L(t)].$$ Recalling the probability density propriety \eqref{tablB}, the above model \eqref{SptialP} satisfies the  mass conservation   hypothesis, i.e. $ \frac{d\rho}{dt}= 0$, as it is required in spatially homogeneous conditions.

\subsection{Well-posedness of the related Cauchy problem }
Let denote
$\mathbb{M}_{Ln}$ a set of matrix endowed with the $1-$norm
$$ \parallel f(t)\parallel _1= \sum_{\ell=1}^L\sum_{i=1}^n  \mid f_{i}(t)\mid , \quad  f=(f_{i}^\ell) \in \mathbb{M}_{Ln}.$$
We introduce the linear space $X_T=C ([0,T]; \mathbb{M}_{Ln} )$ of the continuous functions $f=f(t):[0,T]\to \mathbb{M}_{Ln}$ for $T>0$ equipped 
with the infinity norm
$$\parallel f \parallel _\infty = \displaystyle\sup_{t \in [0,T] } \parallel f(t)\parallel _1.$$ \\ 
Note that $(X_T,\parallel .\parallel _\infty )$ is a real Banach space. 

{Well-posedness of the spatially homogeneous problem means global in time existence and uniqueness of a solution $f = f(t)$ to the Cauchy problem (\ref{SptialP}). These results pass through two steps. Firstly, we prove the local existence and uniqueness in time of the solution $f$ in $X_T$ for a certain $T > 0$. Secondly, we extend the obtained result to a global solution defined for all $t > 0$. }

We consider the following assumptions to prove the well-posedness of the spatially homogeneous problem \eqref{SptialP}:
$$ (\textbf{\textit{H1}})\qquad {\mathcal{A}_{hp,r} ^{i,\ell}\geq0},\qquad \sum_{\ell=1}^L\sum_{i=1}^n \mathcal{A}_{hp,r} ^{i,\ell}=1 \quad \forall r =1, \dots , L\quad \forall h,p =1, \dots , n $$
whenever $0\leq \rho_\ell{\leq}\frac{1}{L}$ and 
$$(\textbf{\textit{H2}}) \qquad \exists C_{\eta} > 0  \quad /  \quad  0<\eta^\ell[\rho_\ell] \leq C_{\eta}^\ell,\quad \hbox{when} \quad  0 \leq \rho_\ell {\leq} \frac{1}{L},\quad \ell=1,\dots,L. $$

\begin{theorem}\label{Therm1}	Assume $(\textbf{H1})$-$(\textbf{H2})$ are satisfied and let $\displaystyle\sum_{\ell=1}^L\sum_{i=1}^n f^\ell_{ij}:=\rho_0\in[0,1]$, then 
	$\exists \hspace{0.2cm} T>0$ such that problem \eqref{SptialP}  admits unique local  non-negative solution $f \in X_T$. Moreover, it satisfies the mass conservation 
	\begin{equation}\label{Mass}
		\parallel f(t)\parallel _1=\rho_0. \end{equation}
\end{theorem}
\noindent We start by giving some estimates on the nonlinear operator $\mathcal{J}$ given by:
$$\mathcal{J}=\mathcal{J}_{i}^\ell:=\sum_{r\in I_c}\eta^r[\rho_r(t)]\sum_{h,p=1}^n\mathcal{A}_{hp,r}^{i,\ell}f_{h}^rf_{p}^r-\eta^\ell[\rho^\ell(t)]f_{i}^\ell(t)\rho_\ell(t).$$
\begin{lem}
	Assume $(\textbf{\textit{H1}})$ and $(\textbf {\textit{H2}})$  are satisfied, then $\exists\; C_1>0$ such that 
	$$ \parallel \mathcal{J}(f)\parallel _1 \leq C_1\parallel f\parallel _1 ^2,$$
	$$ \parallel \mathcal{J}(f)-\mathcal{J}(g)\parallel _1 \leq C_1(\parallel f\parallel _1  +\parallel g\parallel _1)\parallel f-g\parallel _1.$$
\end{lem}

\noindent Next, we define the following operator $$\phi(f)(t)=\int_0 ^t \mathcal{J}(f)(s) ds. $$ We have the following estimates:

\begin{lem} \label{Lema3} Assume $(\textbf{\textit{H1}})$ and $(\textbf{\textit{H2}})$ are satisfied, then $\exists\, C_2(T)>0$ such that 
	$$ \parallel \phi(f)\parallel _\infty \leq C_2\parallel f\parallel _{\infty} ^2,$$
	$$ \parallel \phi(f)-\phi(g)\parallel _\infty \leq C_2(\parallel f\parallel _\infty +\parallel g\parallel _\infty)\parallel f-g\parallel _\infty.$$
\end{lem}
\noindent\textbf{Proof of Theorem 1} First, we rewrite problem \eqref{SptialP} as an integral form 
$$f=\mathcal{N}(f)$$
where the mild equation is as follows
{\small\begin{eqnarray}\label{MS}
	\displaystyle	(\mathcal{N}(f))_{i}^\ell = f_{i,0}^\ell + \int_0 ^t\Big(\sum_{r=1}^L\eta^r[\rho_r]\sum_{h,p=1}^n   \mathcal{A}_{hp,r} ^{i,\ell} f_{h}^r(s)f^r_{p}(s)- f^\ell_{i}(s) \eta^\ell[\rho_{\ell}]\rho_\ell(s)\Big)ds.
\end{eqnarray}}
Proving the uniqueness and existence of solution of the above mild equation in $X_T$ implies to finding a fixed point of $\mathcal{N}$. This is given by the following lemma
\begin{lem}\label{Lema4}
	$\mathcal{N}$ maps the ball $B_{X_T}(0,2\,\parallel f_0\parallel_1)$ into itself.
\end{lem}
\noindent Let $f \in X_T, \quad R>0$ such that $\parallel f\parallel _\infty \leq R$. From Lemma \ref{Lema3} we obtain
$$ \parallel \mathcal{N}(f)\parallel _\infty \leq \parallel f_0\parallel _1 +T C_1 \parallel f\parallel_\infty ^2,$$
$$ \parallel \mathcal{N}(f)\parallel _\infty \leq R +T C_1 R^2.$$
We chose $T$ such that 
$$T C_1 R^2< \frac{R}{4},$$
this implies
$$T < \frac{1}{4C_1  \parallel f_0\parallel _1}.$$
Then,
$$ \parallel \mathcal{N}(f)\parallel _\infty \leq \parallel f_0\parallel _1 + \frac{1}{4C_1  \parallel f_0\parallel _1} C_1 \parallel f\parallel _\infty ^2,$$
$$\parallel \mathcal{N}(f)\parallel_\infty \leq \parallel f_0\parallel_1 + \frac{1}{4 \parallel  f_0\parallel_1}R ^2.$$
where $R$ is a solution of  $\frac{1}{4} R^2 - \parallel f_0\parallel _1 R+ \parallel f_0\parallel _1 ^2=0.$
Then, $R=2\parallel f_0\parallel _1$.
\begin{lem}\label{Lema5}
	$\mathcal{N}$ is strictly contraction on a ball $B_{X_T} \big(0,2\parallel f_0\parallel _1\big)$.
\end{lem}
\noindent Let $f , g \in X_T, \; R>0$ such that $\parallel f\parallel_\infty \leq R$ and $\parallel g\parallel_\infty \leq R$. Thanks to Lemma \ref{Lema3}, we have the following estimate
$$ \parallel \mathcal{N}(f)-\mathcal{N}(g)\parallel _\infty \leq  T C_1\big(\parallel f\parallel _\infty +\parallel g\parallel_\infty\big)\parallel f-g\parallel_\infty.$$
It follows that
$$ \parallel \mathcal{N}(f)-\mathcal{N}(g)\parallel_\infty \leq 2R T C_1\parallel f-g\parallel_\infty.$$
We chose $T$ such that
$$T < \frac{1}{2RC_1}.$$
Then, Lemma \ref{Lema5} is achieved.\\
Now, by summing over $\ell,i,j$ and using Assumption (\textbf{\textit{H1}}), one gets the conservation of mass given by equation \eqref{Mass}. Notice that we have the non-negativity of the solution thanks to Lemma \ref{Lema4} under the non-negativity of the initial condition  $f_{i,0}^\ell>0$. {Indeed,
Eq. (\ref{SptialP}) can be rewritten as follows
{\small\begin{equation}\label{*}
			\displaystyle\frac{d}{dt}f_{i}^\ell(t)+f_{i}^\ell(t)L(f)(t)=G_i(f,f)(t):=\sum_{r\in I_c}\eta^r[\rho_r(t)]\sum_{h,p=1}^n\mathcal{A}_{hp,r}^{i,\ell}f_{h}^rf_{p}^r, 
\end{equation}}
where $L(f)(t)=\eta^\ell[\rho_\ell(t)]\rho_\ell(t)$.
 Let denote $\lambda(t)=\int_{0}^{t}L(f)(s)ds$. 
If $f_i^{\ell}$ is a solution of Eq. \eqref{*}, then
$$\displaystyle\frac{d}{dt}\Big(\exp(\lambda(t))f_{i}^\ell(t)\Big)=\exp(\lambda(t))G_i(f,f)(t),$$
which implies 
\begin{equation}\label{**}
	f_i^\ell(t)=\exp(-\lambda(t))f_{i,0}^\ell+\int_{0}^{t}\exp(\lambda(s))G_i(f,f)(s)ds.
\end{equation}
Consequently, eq. \eqref{**} allows to conclude that, given $f_{i,0}^\ell>0$ and the positivity of the integral function, the solution $f_i^{\ell}(t)$ satisfies non-negativity condition  in its domain of existence. }\\
$$\hskip15.5cm\square$$
The following theorem guaranties the existence of global solution of the spatially homogeneous problem \eqref{SptialP}:
\begin{theorem}
	Under the same hypotheses of Theorem \ref{Therm1}, Problem (\ref{SptialP}) admits a unique non-negative global solution $f \in C(\mathbb{R}^+, \mathbb{M}_{Ln})$ satisfying estimate \eqref{Mass}.
\end{theorem}
\noindent\textbf{Proof.} It suffices to apply the same reasoning developed in the proof of Theorem
\ref{Therm1} on the interval $(T^*, 2T^*]$, taking  $f(T^*)$ as new initial condition. Since $f^\ell_{i}(T) \geq 0$ for all
$\ell=1,\dots,L,\;i = 1,\dots, n$, and $\displaystyle\sum_{\ell=1}^L\sum_{i=1}^n f^\ell_{i} = \rho_0 \in[0, 1],$ we are in the same hypotheses of Theorem \ref{Therm1}. Hence, we conclude the existence and uniqueness of a local in time continuous solution on $[T^*, 2T^*]$ satisfying 
\begin{equation}
	\parallel f(t)\parallel _1=\rho_0,\qquad \forall t\in[T^*, 2T^*]. \end{equation}
Iterating this procedure on all intervals of the form $(kT^*, (k+1)T^*], k \in \mathbb{N},$ we can construct the global solution on real positive axis $\mathbb{R}^+$.

\subsection{Numerical simulations}
We present some numerical simulations obtained from the spatially homogeneous problem \eqref{SptialP} to validate proposed modeling of each terms, namely of the table of games. This can be done according to the existing results in the literature. Concretely, Kerner's fundamental diagrams which relates the flux and the average velocity to the density, and the asymptotic property in time of the vehicles on the lanes starting from different initial distributions.
\subsubsection{Kerner's fundamental diagrams}
\begin{figure}
	\centering
	\subfigure{ \includegraphics[height=1.5in ,width=4in]{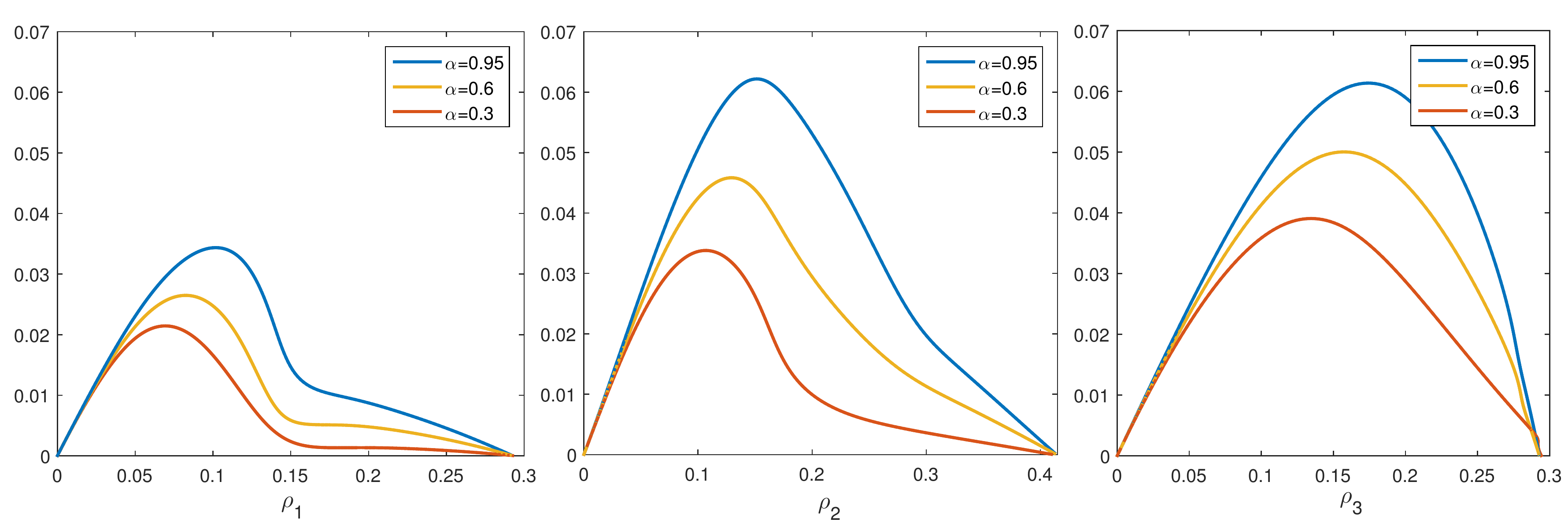}}
	\subfigure{ ~\includegraphics[height=1.5in ,width=4in]{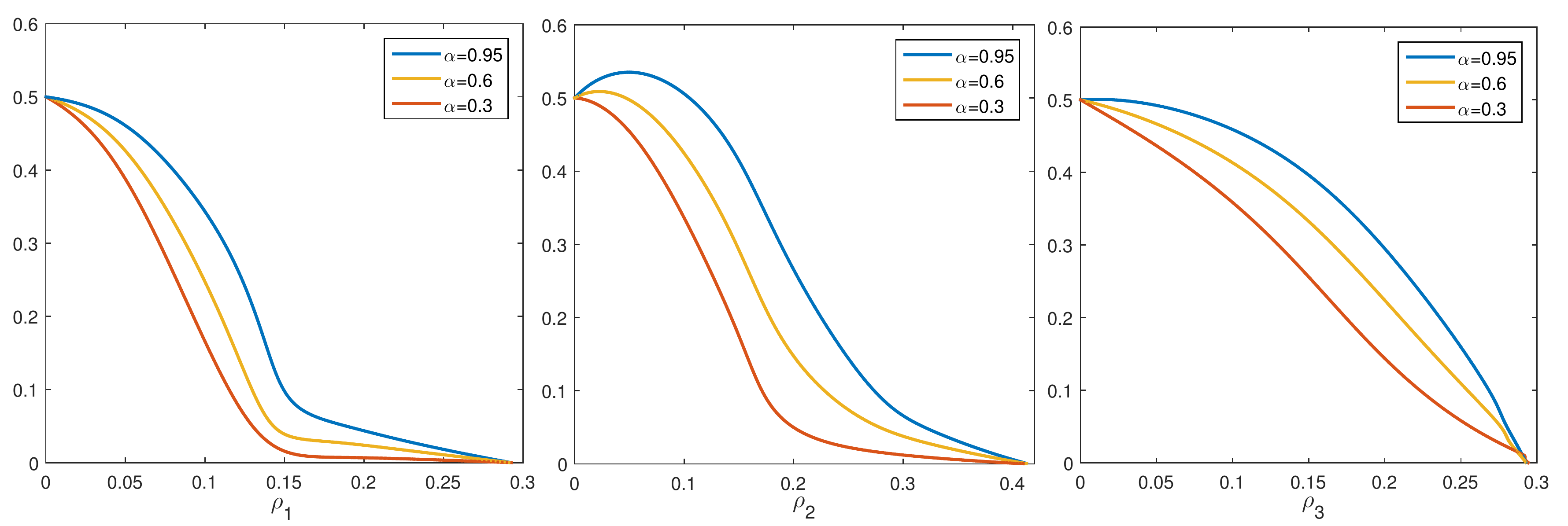}}
	\caption{The first line shows fundamental diagrams: flux $q^\ell$ vs. $\rho_\ell$. The second line shows velocity diagrams: means velocity $U^\ell$ vs. $\rho_\ell$. These numerical results are obtained under various road conditions ($\alpha = 0.95, \alpha = 0.6, \alpha = 0.3$ respectively).}
	\label{5.2}
\end{figure}

\begin{figure}
	\centering
	\subfigure{ \includegraphics[height=1.5in ,width=4in]{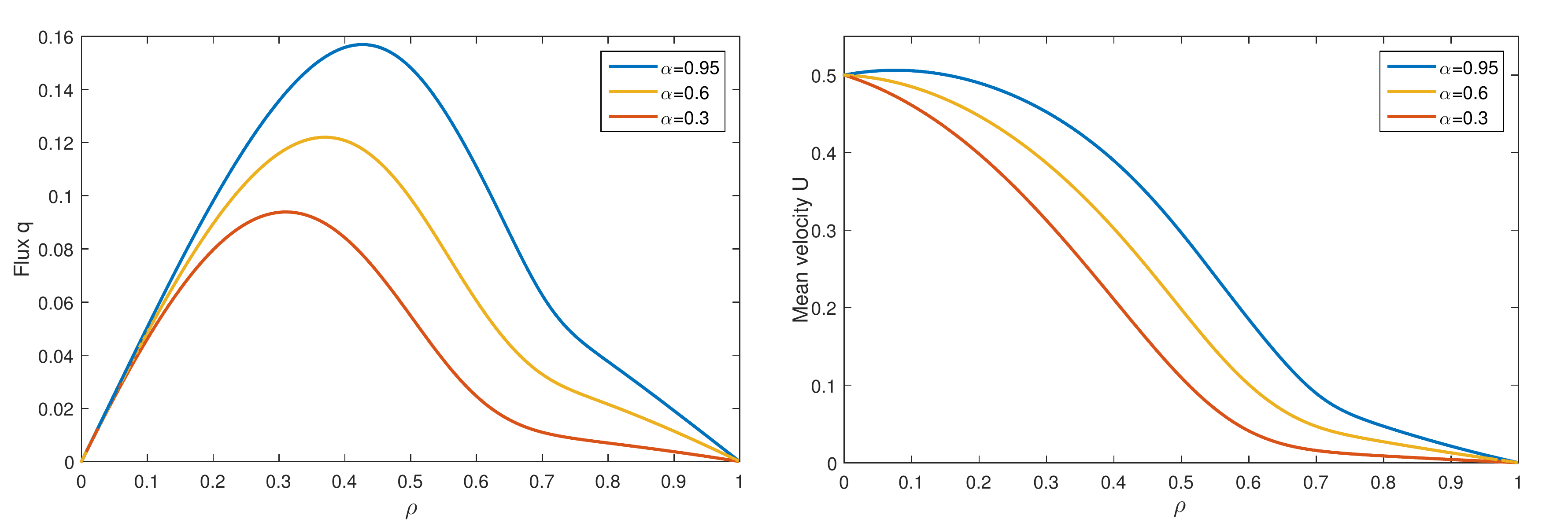}}
	\caption{Global fundamental diagrams obtained under various road conditions ($\alpha = 0.95, \alpha = 0.6, \alpha = 0.3$ respectively): (Left) flux $q$ vs. density $\rho$; (right) means velocity $U$ vs. density $\rho$ .}
	\label{5.3}
\end{figure}

Numerical simulations of Eq.(\ref{SptialP}) have been performed to obtain the fundamental diagrams relating to the average velocity and the macroscopic flux to the density of the vehicles. The problem under consideration has been numerically solved by fixing six-velocity classes $\{v_i\}_{i=1}^6$ and three lanes $L=3$. The result is shown in Figures \ref{5.2}-\ref{5.3} corresponding to three different values of the phenomenological parameter $\alpha$, that is to say $\alpha=0.95, \; \alpha=0.6, \; \alpha=0.3$
and by fixing $\rho^0_{\ell}\in \big[0,\frac{1}{L}\big[, \ell=1,2,3,$ at the initial time.\\
Figure \ref{5.2} shows, respectively, the result of fundamental diagrams for three-lane road, namely the slowest lane, the medium lane, and the fastest lane, and Figure \ref{5.3} shows the global ones. The obtained results can be summarized as follows:
\begin{itemize}
	\item [\textbullet] For low density the flux exhibits an almost linear behavior, which is in agreement with experimental observation reported by Kerner under free-flow conditions;
	\item [\textbullet] The flux becomes markedly nonlinear when the density increases;
	\item [\textbullet] The average speed is initially almost close to the maximum possible, then it drops steeply to zero when the density enters the congested flow range;
	\item[\textbullet] The free flow phase reduces as the environmental conditions worsen;
	\item[\textbullet] The free flow  increases as the vehicles are positioned in the fastest lane.  
\end{itemize}

Figures \ref{5.2}-\ref{5.3} demonstrate that our model is able to catch qualitatively the well-known phase transition from free to congested traffic flow.

\subsubsection{Asymptotic property in time}
We aim to show the asymptotic property in time of the vehicles on the three lanes ($L=3$) beginning from different initial conditions as follows
\begin{equation}
	\begin{array}{llll}
		\rho=0.2:&\rho_{1}\mid_{t=0}=0.2,& \rho_{2}\mid_{t=0}=0,&\rho_{3}\mid_{t=0}=0;\\
		\rho=0.4:&\rho_{1}\mid_{t=0}=0.2,& \rho_{2}\mid_{t=0}=0.2,&\rho_{3}\mid_{t=0}=0;\\
		\rho=0.6:&\rho_{1}\mid_{t=0}=0.3,& \rho_{2}\mid_{t=0}=0.15,&\rho_{3}\mid_{t=0}=0.15;\\
	\end{array}
\end{equation} 

\begin{figure}
	\centering
	\subfigure{ \includegraphics[height=1.5in ,width=4in]{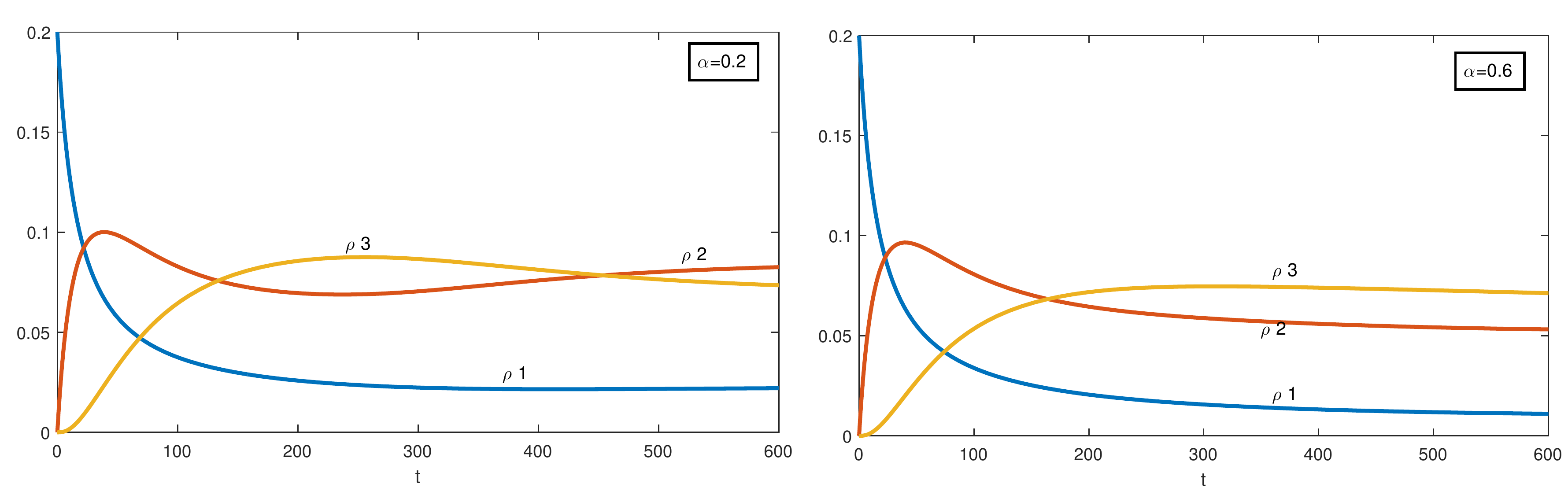}}
	\caption{Global density $\rho=0.2$ at different values of $\alpha$.}
	\label{5.4}
\end{figure}

\begin{figure}
	\centering
	\subfigure{ \includegraphics[height=1.5in ,width=4in]{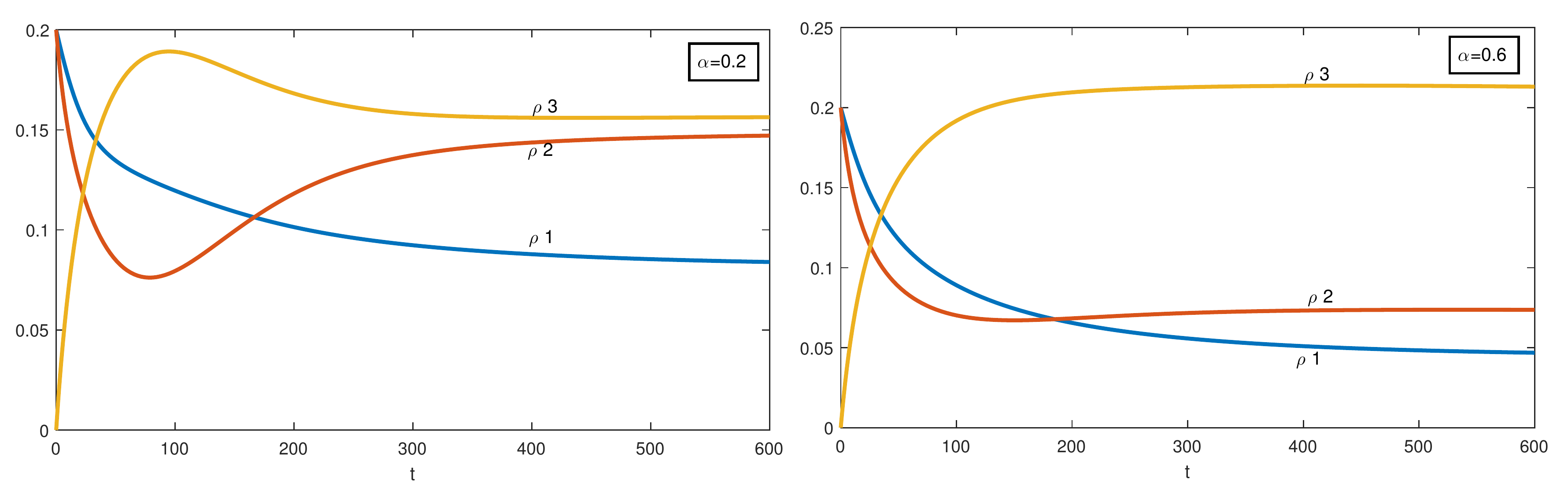}}
	\caption{Global density $\rho=0.4$ at different values of $\alpha$.}
	\label{5.5}
\end{figure}

Figure \ref{5.4} shows the obtained numerical simulations in the case of a low global density of $\rho=0.2$ where the micro-systems started from the slowest lane. We notice that when the equilibrium is achieved, micro-systems are distributed in all lanes. Specifically, the medium lane is the most occupied in the case of bad road conditions $\alpha=0.2$ while the fasted lane is the busiest in the case of good road conditions $\alpha=0.6$. Next, we increase the global density $\rho=0.4$ and we assume the same initial distributions for the first and the second lanes while the third lane is empty. Figure \ref{5.5} shows the obtained results. We notice that the fastest lane is the most occupied lane at equilibrium and when the road conditions are bad the fastest and the medium lanes are occupied all the most in the same way. Finally, we increase the global density $\rho=0.6$ and we assumed that the slowest lane is the busiest. We notice that the fastest lane is the most occupied at equilibrium, see Figure \ref{5.5}. Precisely, we have to following results
\begin{itemize}
	\item Increasing the density increases the rapidity by which the steady density is reached;
	\item Micro-systems tend to occupy the fastest lane. This trend is all the more remarkable as the quality of the road increases and the density is low. On the contrary, we note that when the density is important, the micro-systems tend to occupy the entire line and the state of the road parameter is less relevant;
	\item The same road conditions and the same overall density show a tendency towards a balance insensitive to the initial distribution on the lanes.	
\end{itemize}
\begin{figure}
	\centering
	\subfigure{ \includegraphics[height=1.5in ,width=4in]{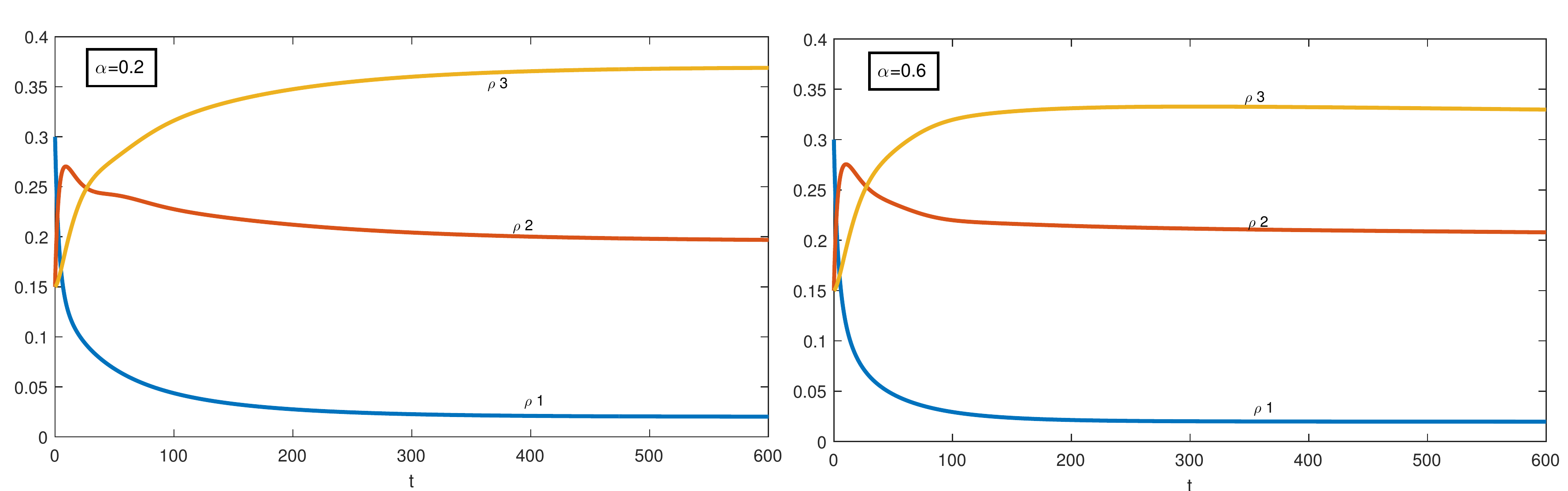}}
	\caption{Global density $\rho=0.6$ at different values of $\alpha$.}
	\label{5.6}
\end{figure}

\newpage
\section{The spatially inhomogeneous problem}\label{Sect5}
This section is devoted to present the numerical simulations of the spatially inhomogeneous problem, which describes the spatial and temporal evolution of the traffic flow. Precisely, we aim to validate our proposed modeling by showing the emerging clusters. Herein, we neglect the term of the external force $\mathcal{T}_{i}^\ell$, and we consider the periodic conditions. Thus, the final form of model (\ref{Struct}) is as follows
{\small\begin{equation}\label{InhomgP}
	\left\{\begin{array}{ll}
		\displaystyle\frac{\partial f_{i}^\ell}{\partial t}+v_i\;\partial_xf_{i}^\ell=\mathcal{J}_{i}^\ell[\textbf{f},\textbf{f}](t,x),&\forall (t,x)\in [0,T]\times]0,1[,\\\\
		f_{i}^\ell(0,x)=f_{i}^{\ell 0}(x), & \forall x\in[0,1],\\\\
		f_{i}^\ell(t,0)=f_{i}^\ell(t,1), &\forall t\in [0,T],
	\end{array}\right.
\end{equation}}
for $\ell=1,\dots,L$ and $i=1,\dots,n$,  
where
{\small\begin{equation*}
	\begin{array}{l}
		\displaystyle\mathcal{J}_{i}^\ell[\textbf{f},\textbf{f}](t,x)=\sum_{r\in I_c}\sum_{h,p=1}^{n}\int_{x}^{x+\xi}\eta^r[\rho^\star_r(t,x^*)]\mathcal{A}_{hp,r}^{i,\ell}[\rho_{1}^\star,\dots,\rho_L^\star;\alpha](t,x^*)\\
		\hskip 3cm \times f_{h}^r(t,x)f_{p}^r(t,x^*)\omega^r(x,x^*)dx^*
		\\
		\hskip 3cm \displaystyle-f_{i}^\ell(t,x)\int_{x}^{x+\xi}\eta^\ell[\rho^\star_\ell(t,x^*)]\rho_\ell(t,x^*)\omega^\ell(x,x^*)dx^*.
	\end{array}
\end{equation*}}

\begin{figure}
	\centering
	\subfigure[]{ \includegraphics[height=1.5in ,width=2in]{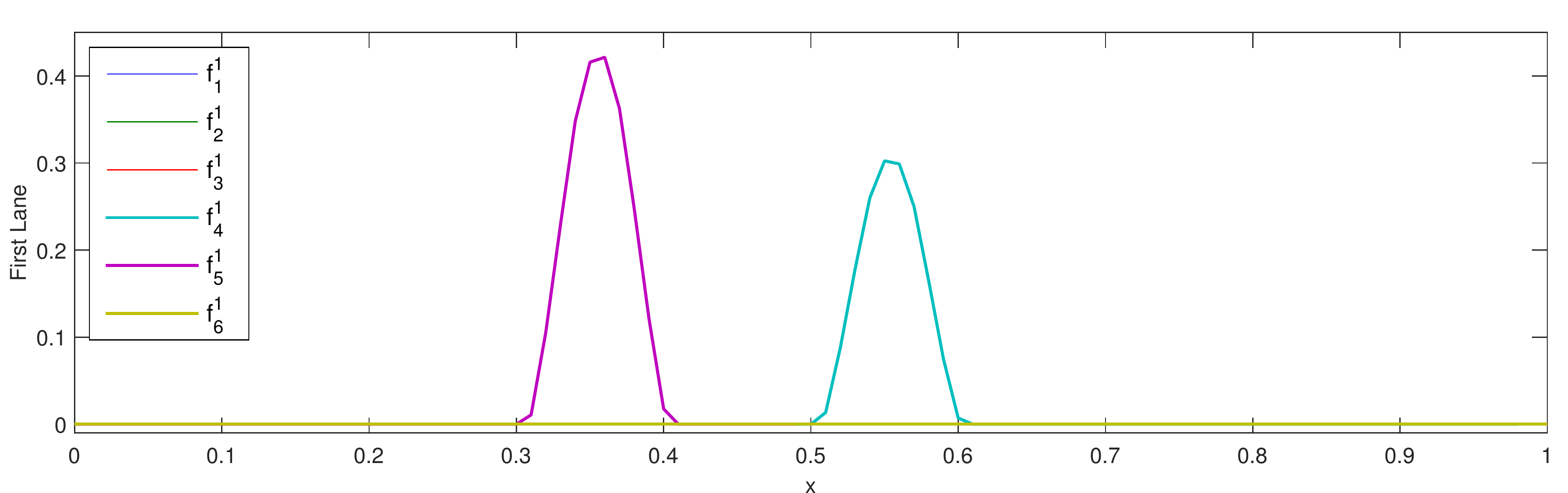}}
	~\subfigure[]{ \includegraphics[height=1.5in ,width=2in]{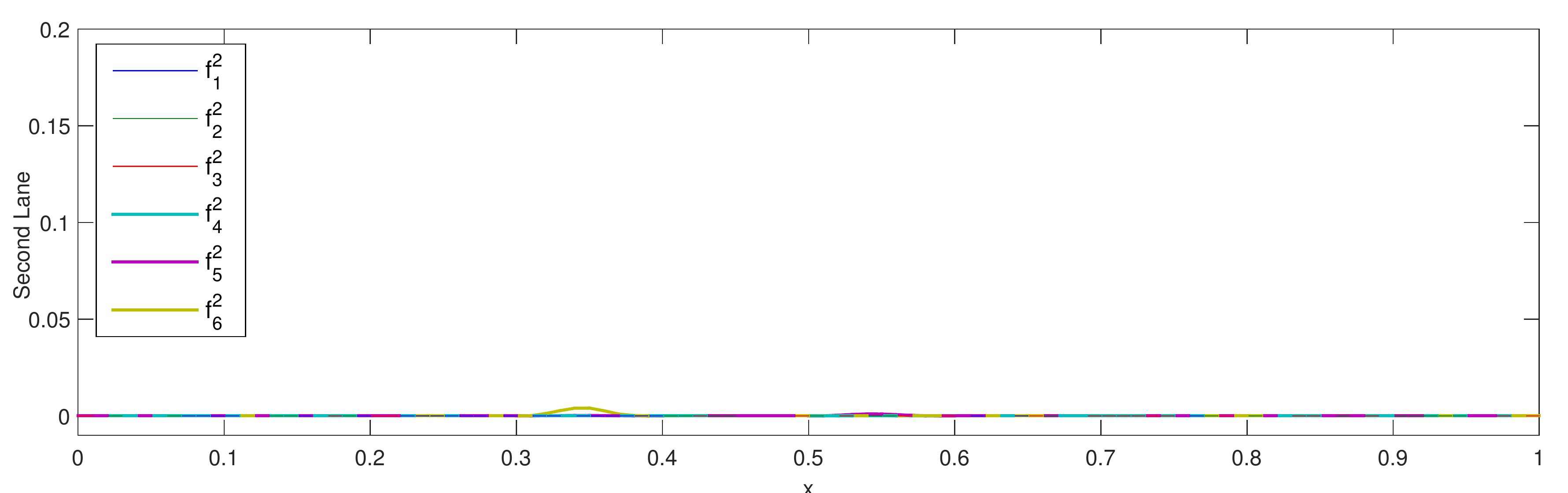}}\\
	{Snapshot of $f_{i}^1$ (left) and $f_i^2$ (right) for $i=1,\dots,6$ at time $t=0.01$.}\\	
	\subfigure[]{ \includegraphics[height=1.5in ,width=2in]{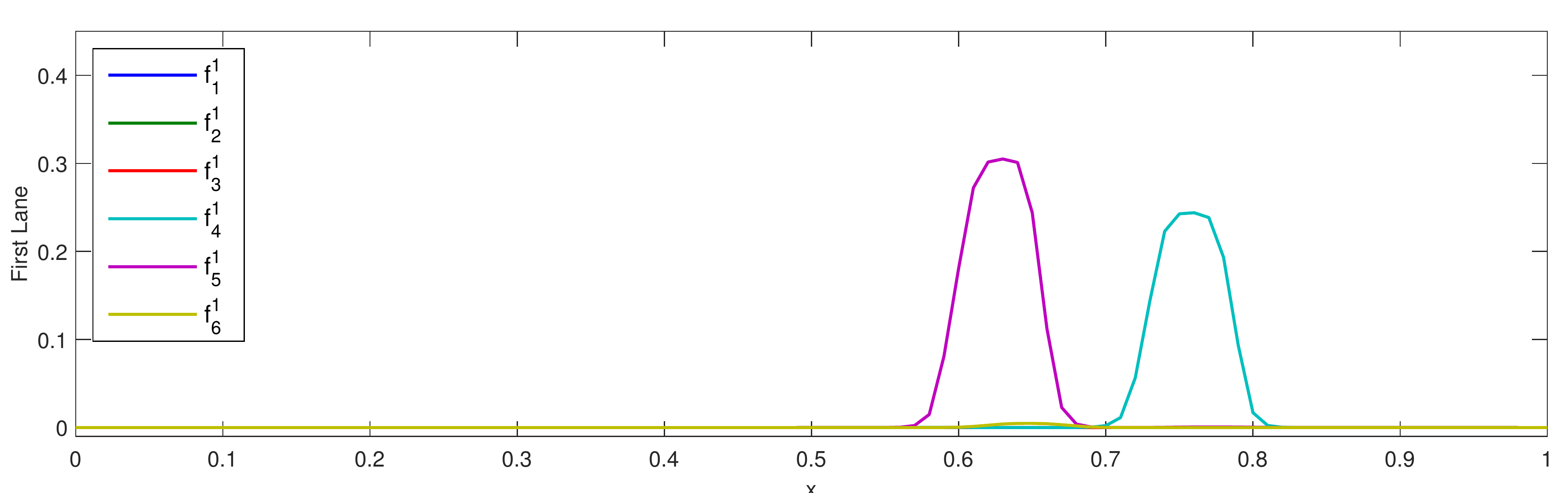}}
	~\subfigure[]{ \includegraphics[height=1.5in ,width=2in]{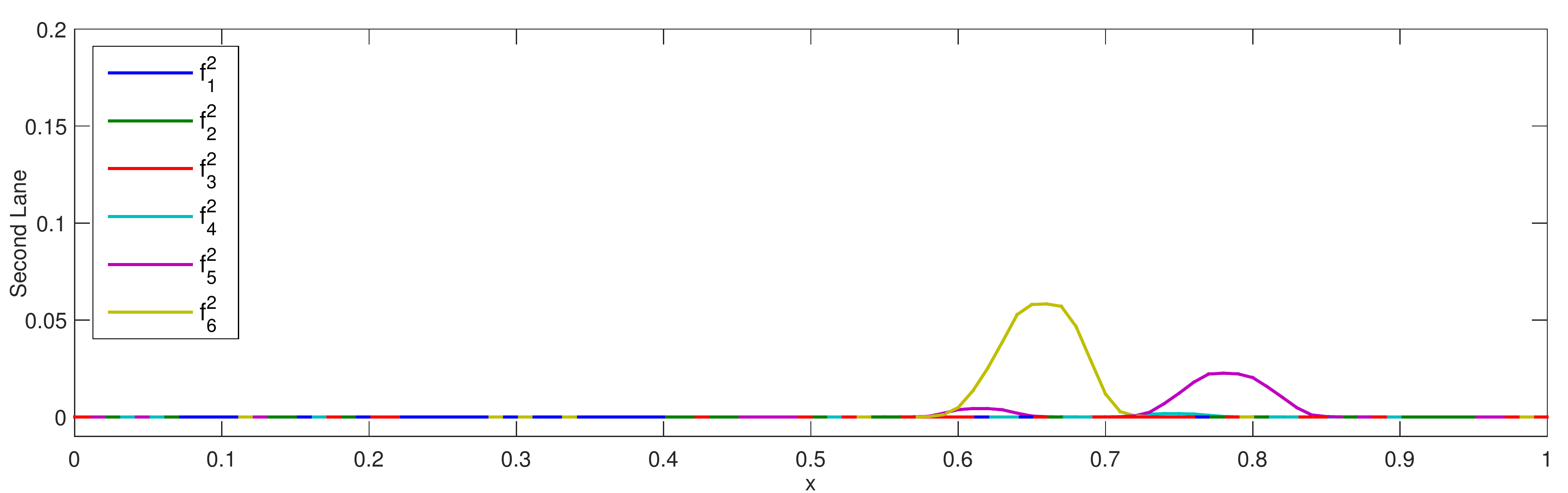}}\\
	{Snapshot of $f_{i}^1$ (left) and $f_i^2$ (right) for $i=1,\dots,6$ at at time $t=0.35$.}\\
	\subfigure[]{ \includegraphics[height=1.5in ,width=2in]{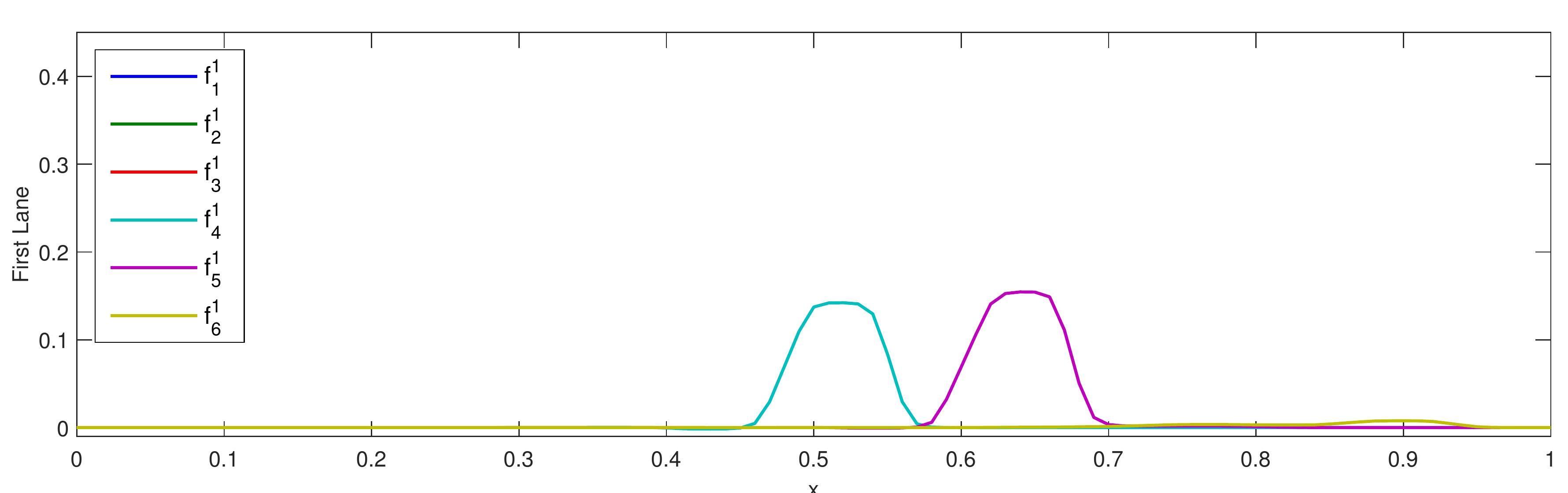}}
	~\subfigure[]{\includegraphics[height=1.5in ,width=2in]{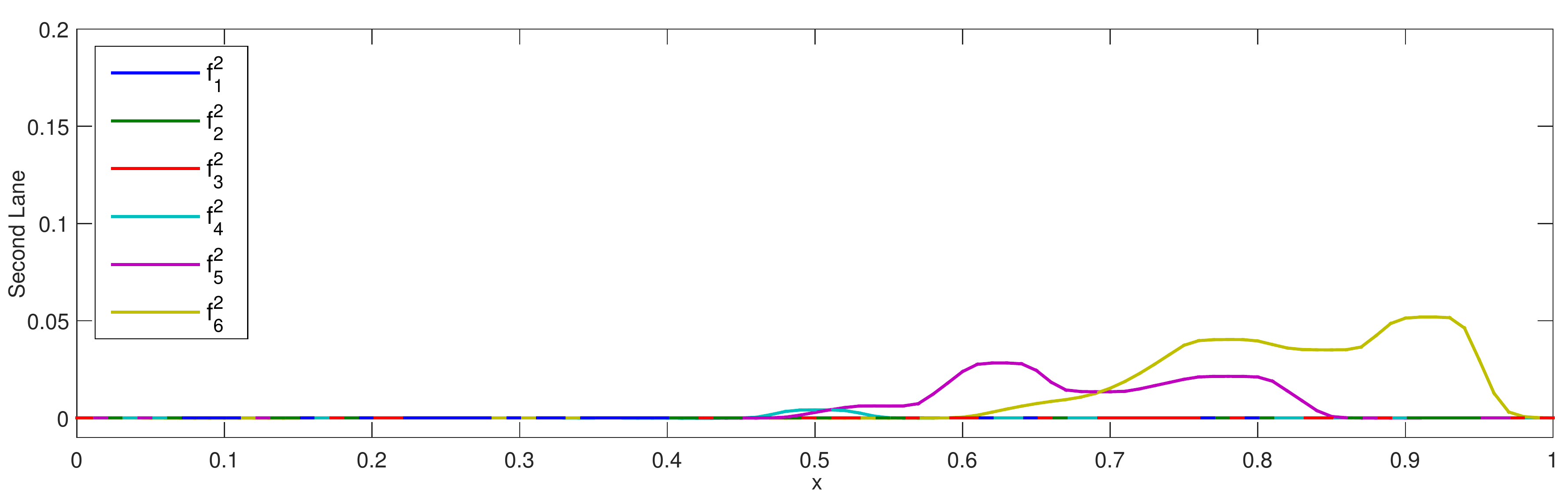}}\\
	{Snapshot of $f_{i}^1$ (left) and $f_i^2$ (right) for $i=1,\dots,6$ at at time $t=1.61$.}\\
	\caption{Evolution of the micro-system in the slowest $f_i^1$ and the fastest lanes $f_i^2$ for $i=1,\cdots,6$ under good road conditions $(\alpha=0.95)$ at successive time $t=0.001,\,0.35,\, 1.61$.}
	\label{1.5}
\end{figure}

We show the numerical simulations of the flow dynamics in two lanes: In the slowest lane, we consider two clusters of vehicles traveling with different speeds while the fastest lane is empty. More precisely, we take
$$f_5^1(t,x)=70\sin^2\big(10\pi(x-0.3)(x-0.4)\big), \qquad x\in[0.3,0.4],$$
$$f_4^1(t,x)=50\sin^2\big(10\pi(x-0.5)(x-0.6)\big), \qquad x\in[0.5,0.6].$$
The hyperbolic system (\ref{InhomgP}) has been integrated by the finite volume method based on the slop limiters technique \cite{[leV]}. The simulations were developed, for different values of the road conditions  $\alpha$, by fixing six velocity classes. Namely, Figure \ref{1.5} refers to a high quality road environment $\alpha=0.95$, we notice that:

\noindent \textbullet $\;$ In the slowest lane, the fast group of vehicles after having approached the group of the slow ones, as indicated in Figure \ref{1.5}-(c), leave them behind, while a small group of vehicles with the velocity $v_n$ is taking the lead, see Figure \ref{1.5}-(e);\\
\noindent \textbullet $\;$ In the fastest lane,  groups of vehicles appeared. The groups of vehicles in the slowest lane have increased their speed when they took the fastest lane (Figure \ref{1.5}-(d)). The group with a greater speed has taken the lead over the other group. The vehicles in the two groups are no longer as united as they used to be in the slowest lane (Figure \ref{1.5}-(f)), and that is due to the fact that some of the vehicles now have the chance to accelerate.\\ 

\begin{figure}
	\centering
	\subfigure[]{\includegraphics[height=1.5in ,width=2in]{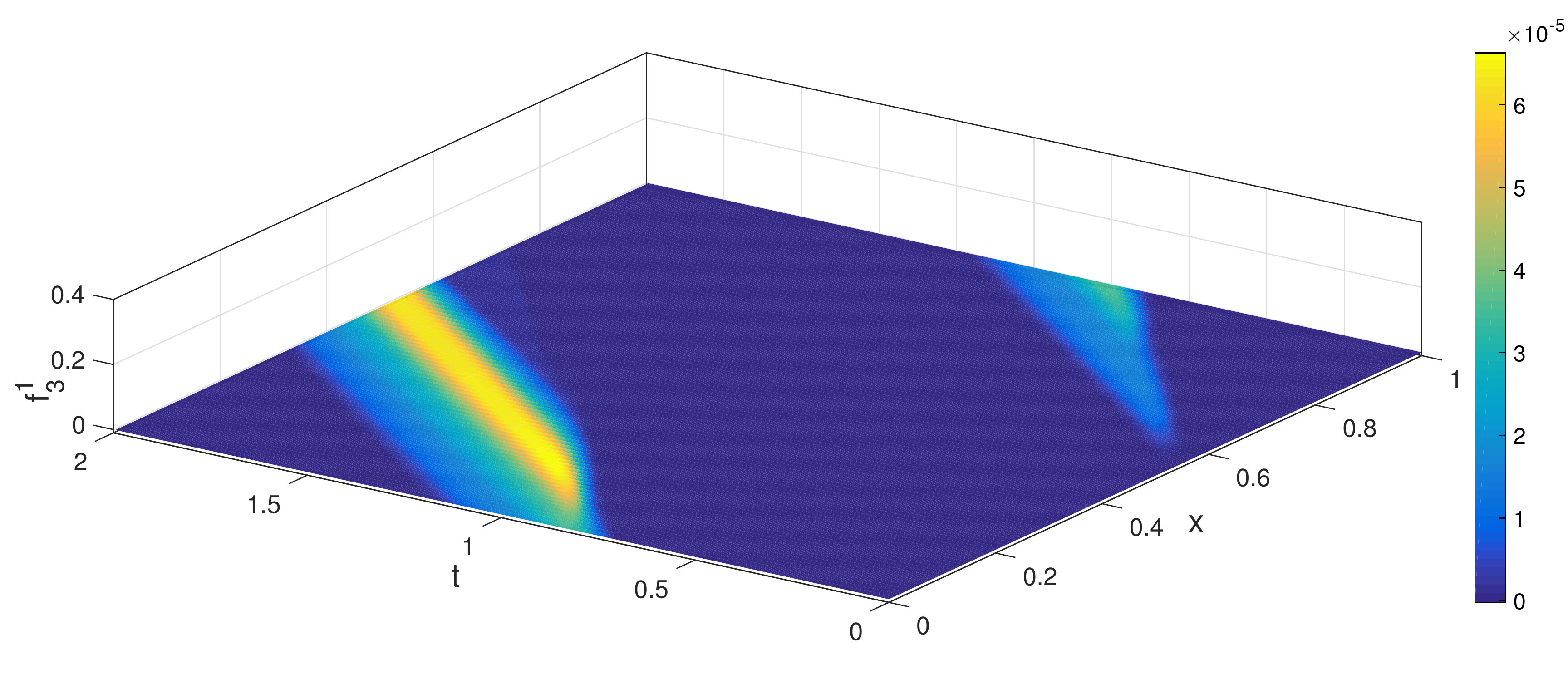}}
	~\subfigure[]{\includegraphics[height=1.5in ,width=2in]{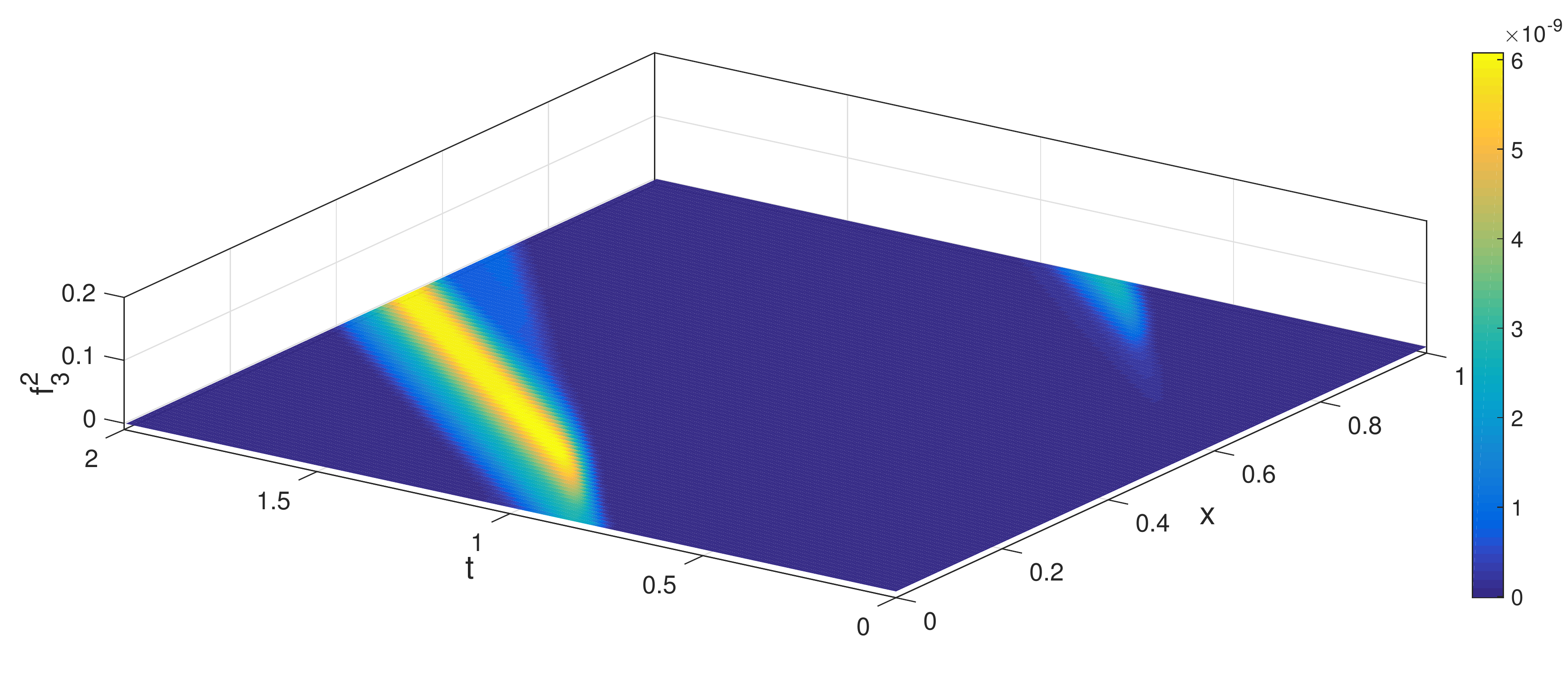}}\\
	{Evolution of $f_{3}^1$ (left) and $f_3^2$ (right).}\\
	\subfigure[]{ \includegraphics[height=1.5in ,width=2in]{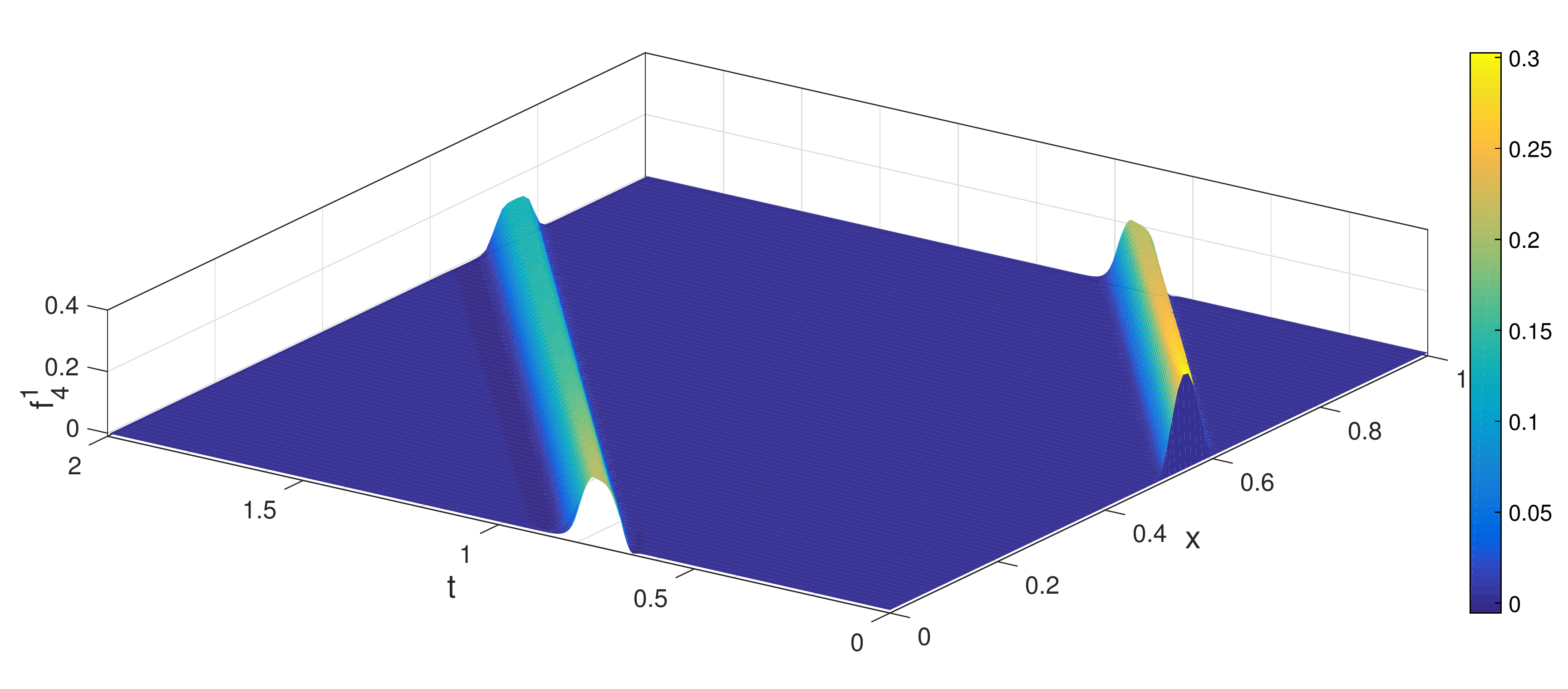}}
	~\subfigure[]{ \includegraphics[height=1.5in ,width=2in]{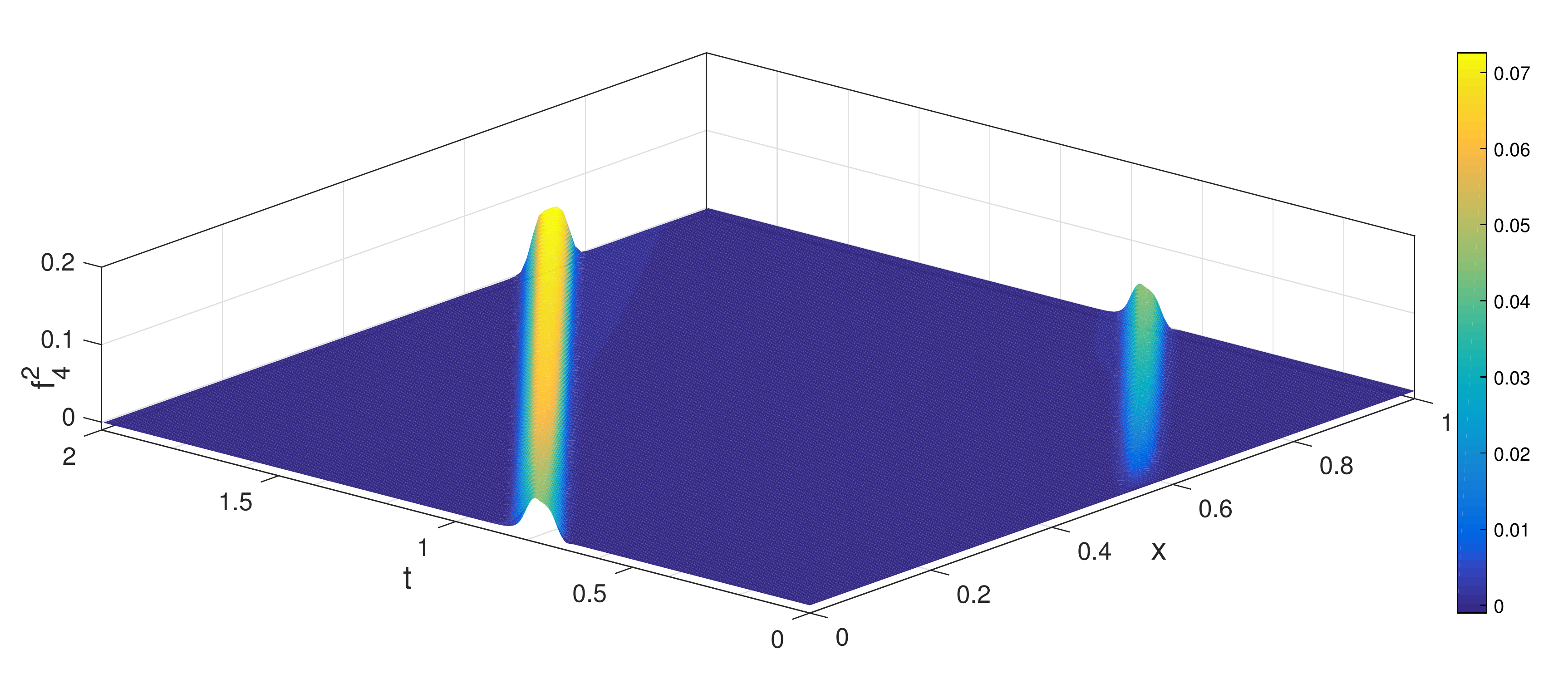}}\\
	{Evolution of $f_{4}^1$ (left) and $f_4^2$ (right).}\\
	\subfigure[]{\includegraphics[height=1.5in ,width=2in]{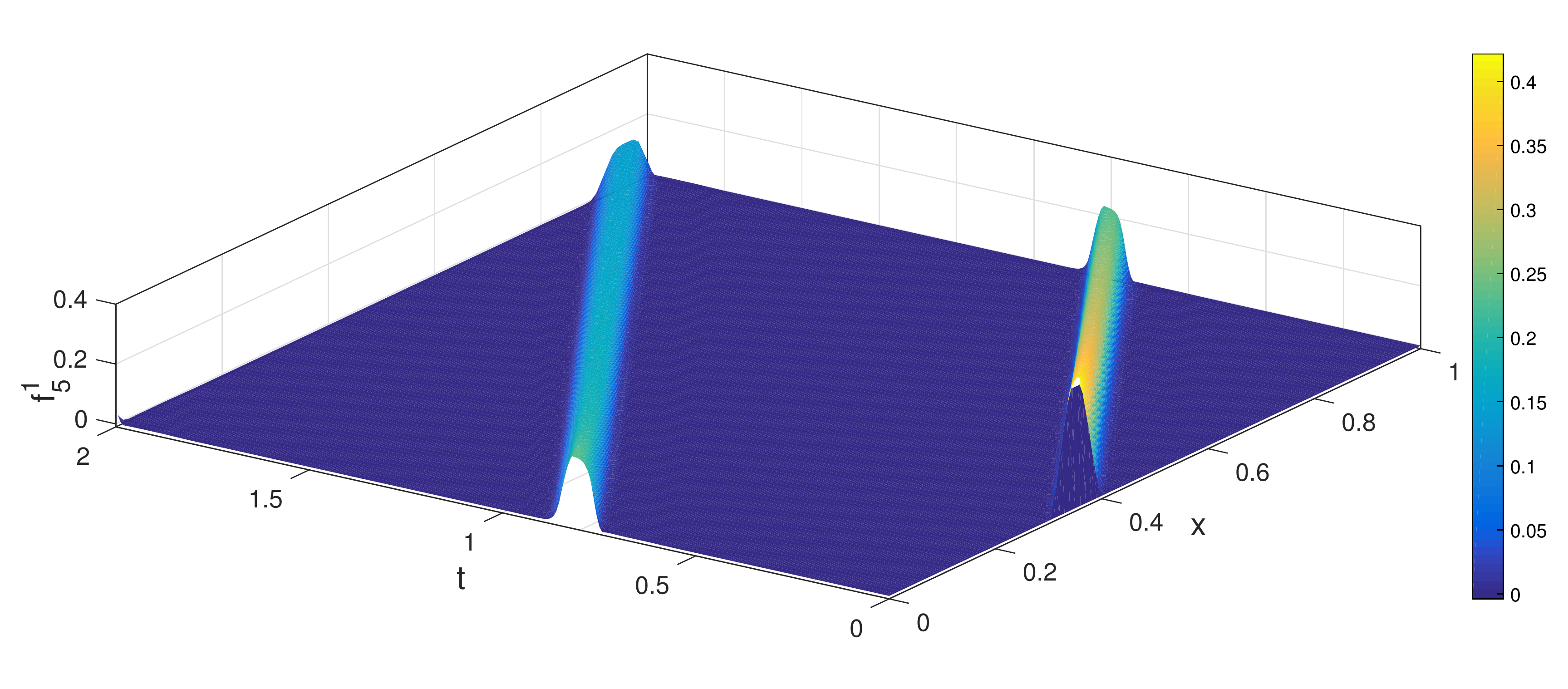}}
	~\subfigure[]{\includegraphics[height=1.5in ,width=2in]{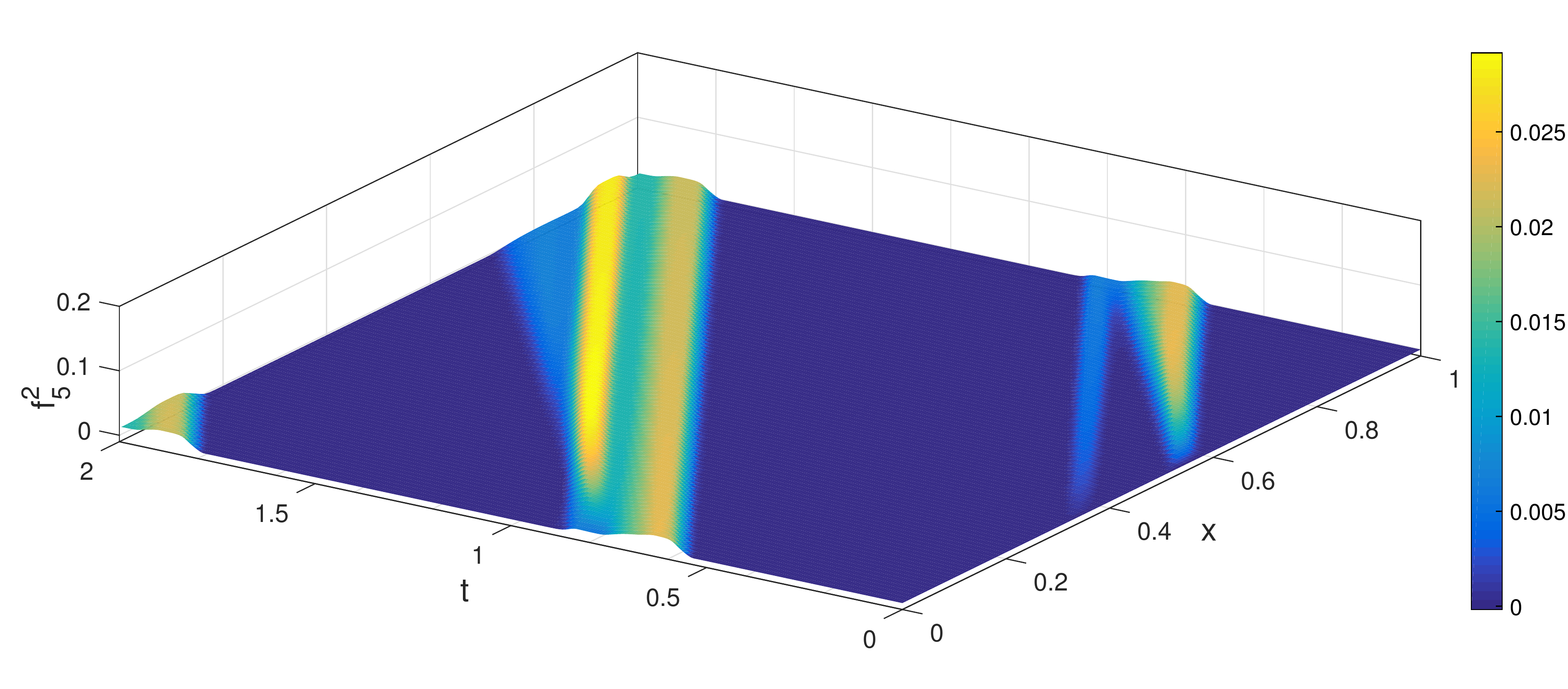}}\\
	{Evolution of $f_{5}^1$ (left) and $f_5^2$ (right).}\\
	\subfigure[]{\includegraphics[height=1.5in ,width=2in]{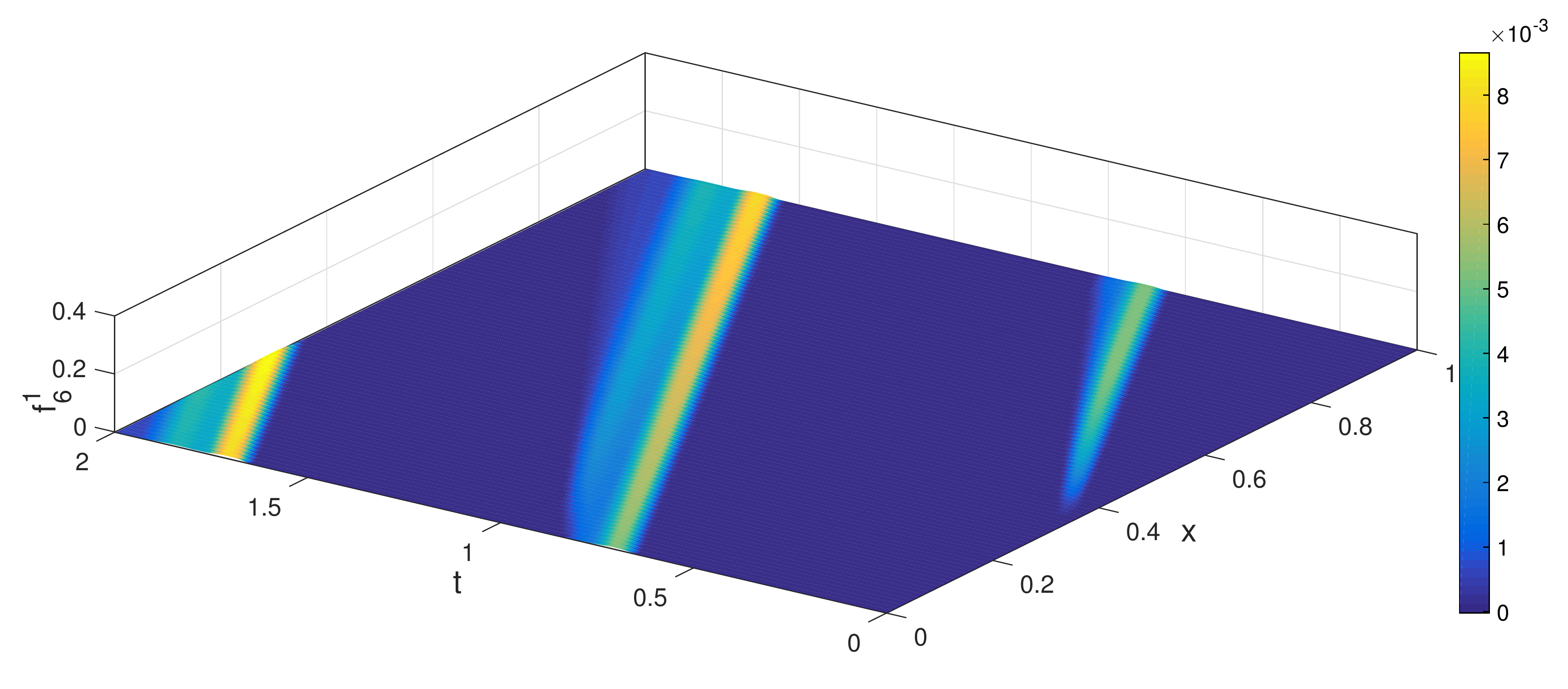}}
	~\subfigure[]{\includegraphics[height=1.5in ,width=2in]{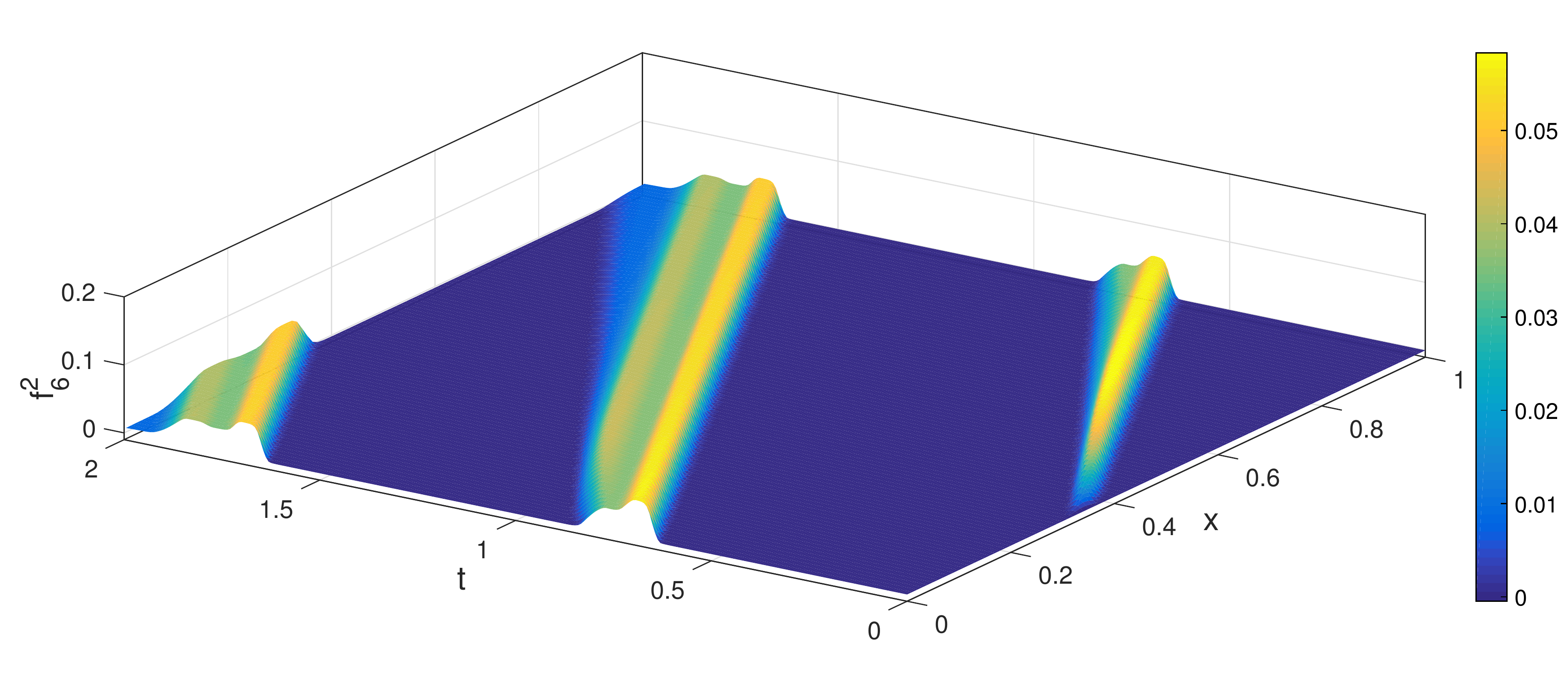}}
	\\{Evolution of $f_{6}^1$ (left) and $f_6^2$ (right).}
	\caption{Evolution of the micro-system in the in slowest lane (left) and in the fastest lane (right) under good road conditions $\alpha=0.95$.}
	\label{1.6}
\end{figure}

Figure \ref{1.6} shows the evolution time $t$ of the micro-systems $f_i^1$ (left) and $f_i^2$ for $i=3,\dots,6$ under the good road conditions ($\alpha=0.95$). We observe that the number of vehicles with velocities $v_{n-2}$ and $v_{n-1}$ have decreased over time as indicated in Figure \ref{1.6}-(c) and \ref{1.6}-(e). Moreover, a group of negligible vehicles has decelerated. Whereas, another group of a relating greater number of vehicles has accelerated (see Figure \ref{1.6}-(a) and \ref{1.6}-(g)). On the other hand, Figure \ref{1.6} shows the number of vehicles of each group that has taken the fastest lane. We notice that the groups with velocities $v_{n-1}$ and $v_n$ (see Figure \ref{1.6}-(f) and \ref{1.6}-(h)) are the only groups that are greatly increasing in number with the time being. Moreover, we observe that the vehicles of groups are no longer in the united groups as mentioned before.

\begin{figure}
	\centering
	\subfigure[]{\includegraphics[height=1.5in,width=2in]{L1T1.pdf}}~\subfigure[]{\includegraphics[height=1.5in ,width=2in]{L2T1.pdf}}
	\\{Snapshot of $f_{i}^1$ (left) and $f_i^2$ (right) for $i=1,\dots,6$ at time $t=0.01$.}\\	
	\subfigure[]{\includegraphics[height=1.5in ,width=2in]{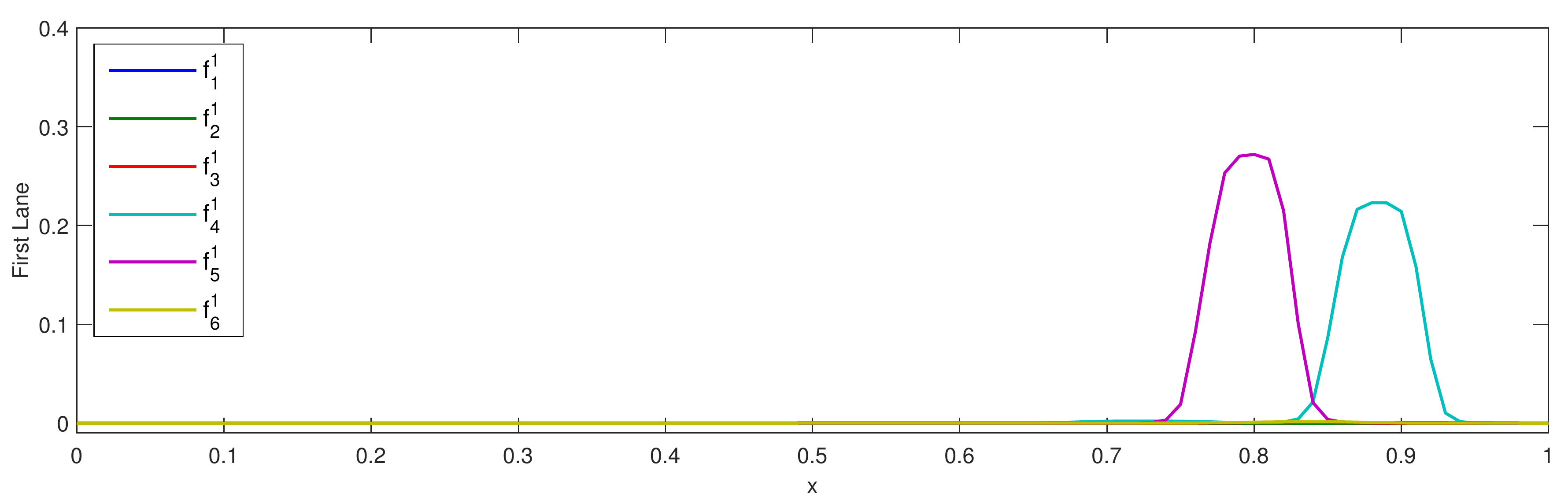}}
	~\subfigure[]{\includegraphics[height=1.5in ,width=2in]{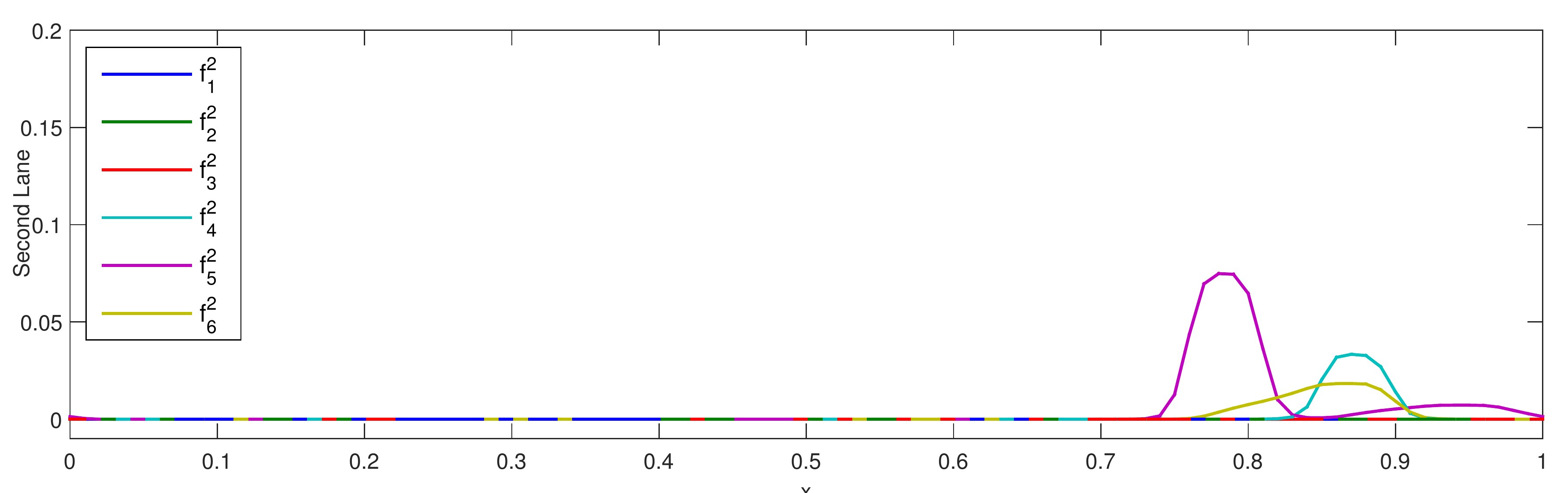}}
	\\{Snapshot of $f_{i}^1$ (left) and $f_i^2$ (right) for $i=1,\dots,6$ at at time $t=0.35$.}\\
	\subfigure[]{\includegraphics[height=1.5in ,width=2in]{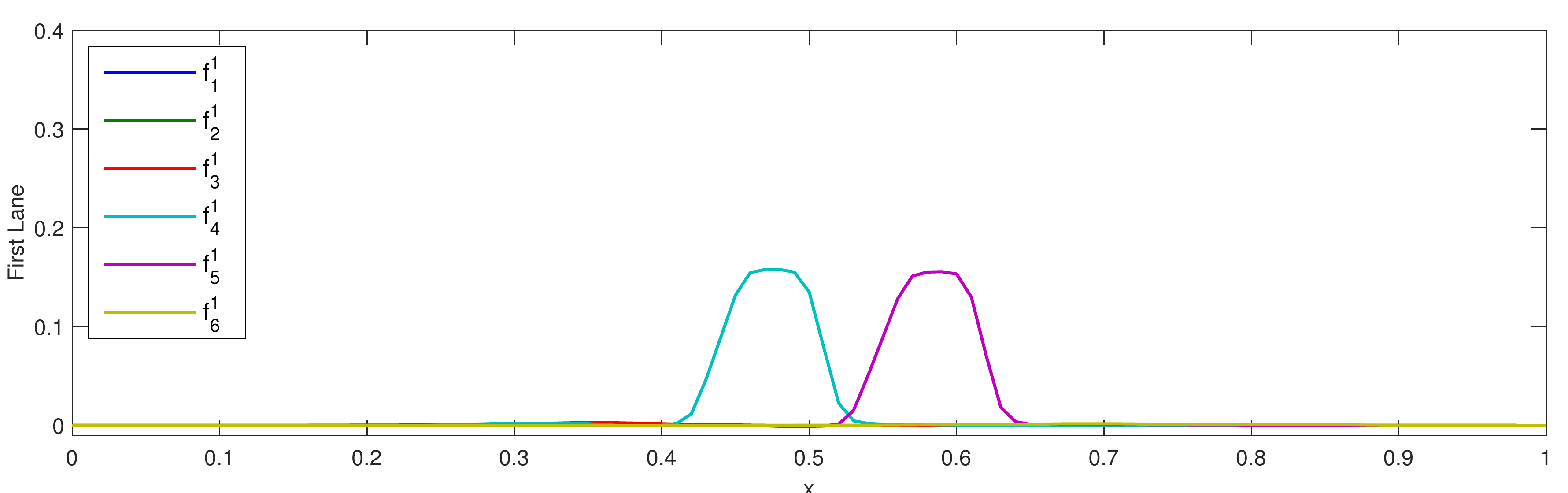}}
	~\subfigure[]{ \includegraphics[height=1.5in ,width=2in]{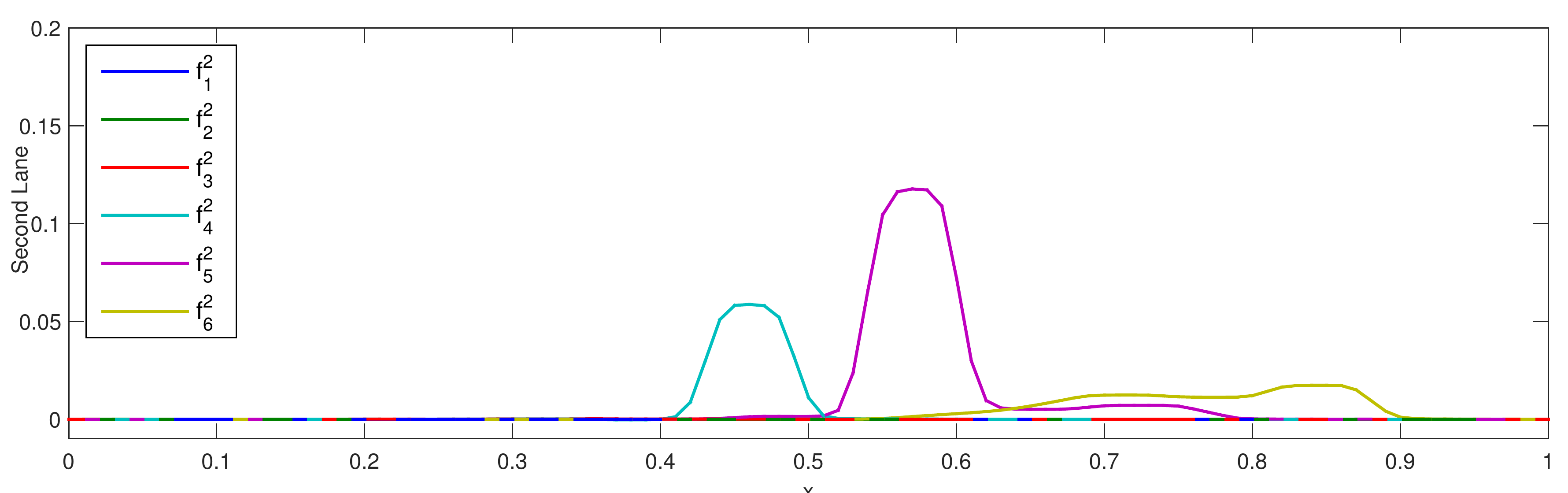}}
	\\{Snapshot of $f_{i}^1$ (left) and $f_i^2$ (right) for $i=1,\dots,6$ at at time $t=1.61$.}\\
	\caption{Evolution of the micro-system in the slowest $f_i^1$ and the fastest lanes $f_i^2$ for $i=1,\cdots,6$ under bad road conditions $(\alpha=0.3)$ at successive time $t=0.001,\,0.35,\, 1.61$.}
	\label{1.7}
\end{figure}

\begin{figure}
	\centering
	\subfigure[]{\includegraphics[height=1.5in ,width=2in]{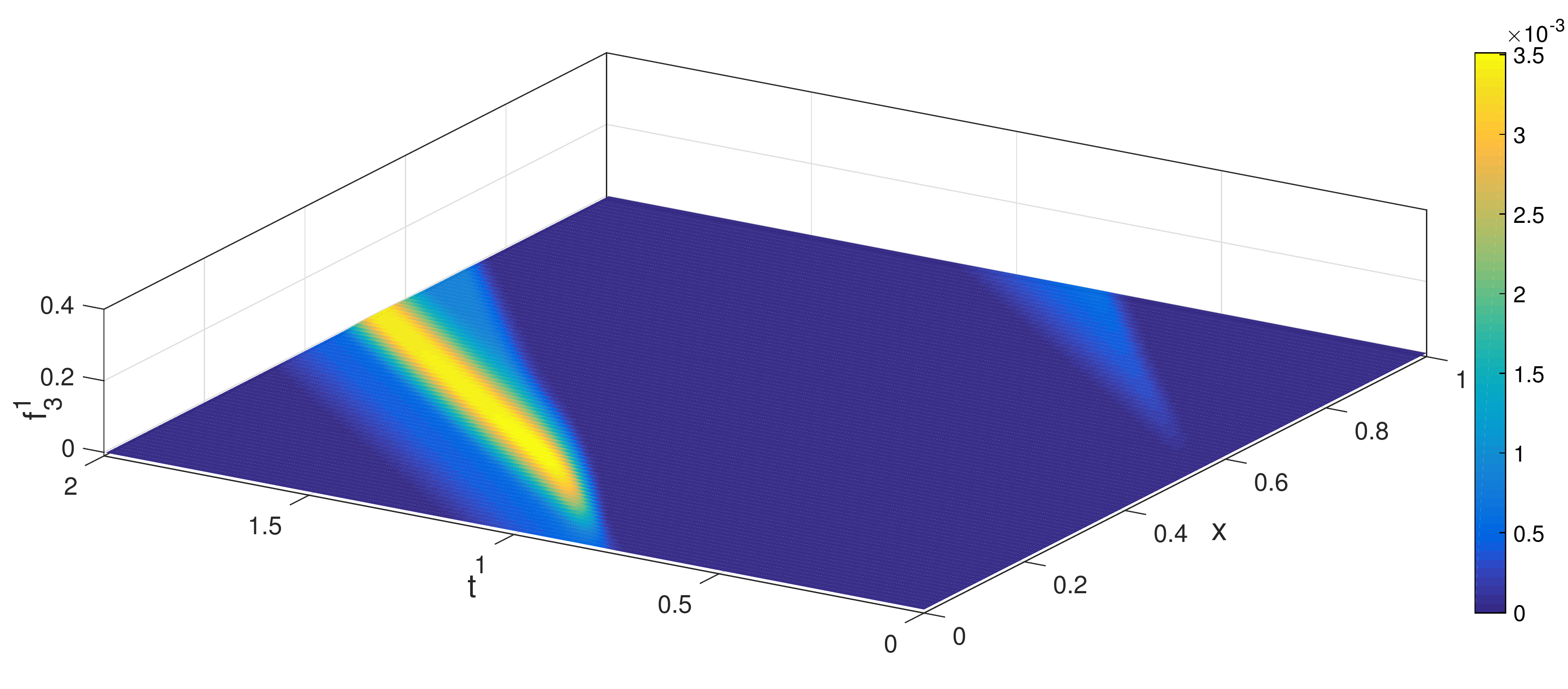}}
	~\subfigure[]{\includegraphics[height=1.5in ,width=2in]{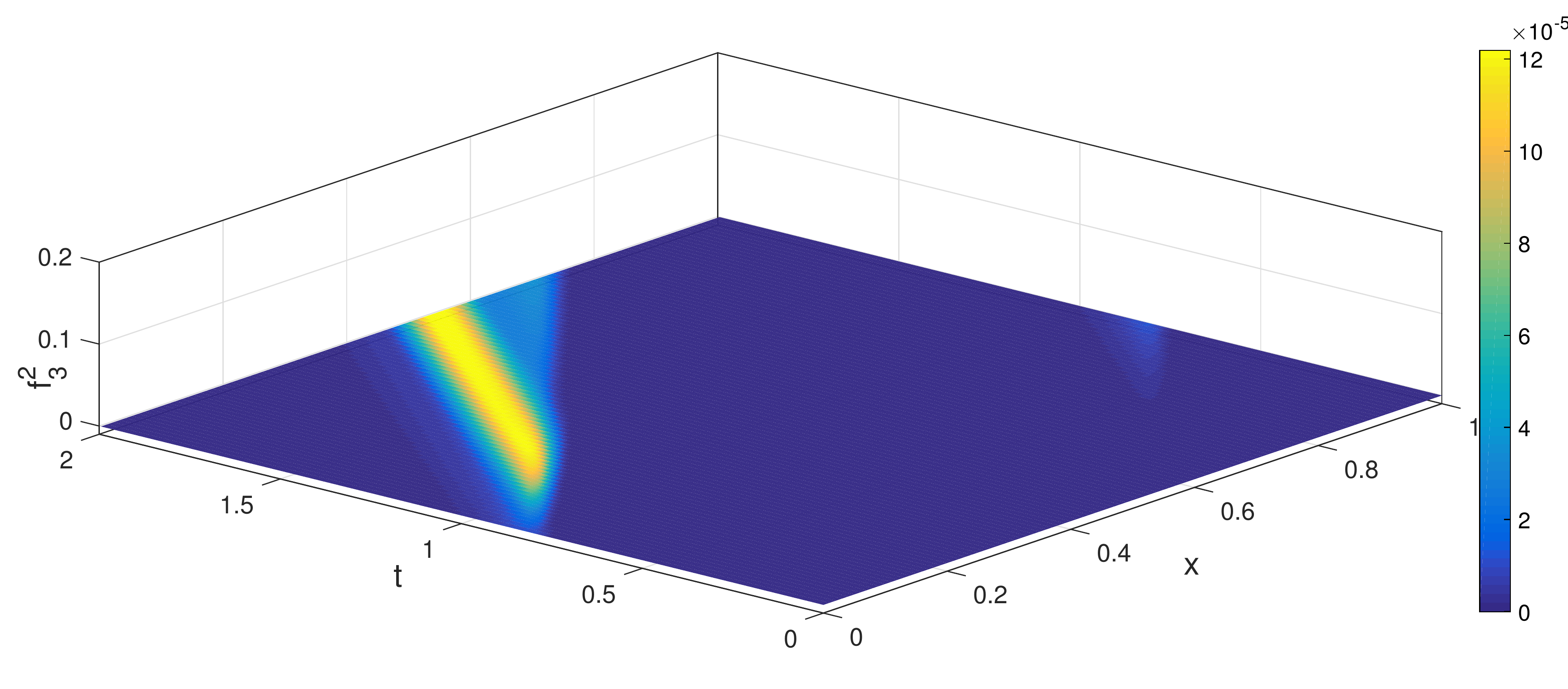}}
	\\{Evolution of $f_{3}^1$ (left) and $f_3^2$ (right).}\\
	\subfigure[]{\includegraphics[height=1.5in ,width=2in]{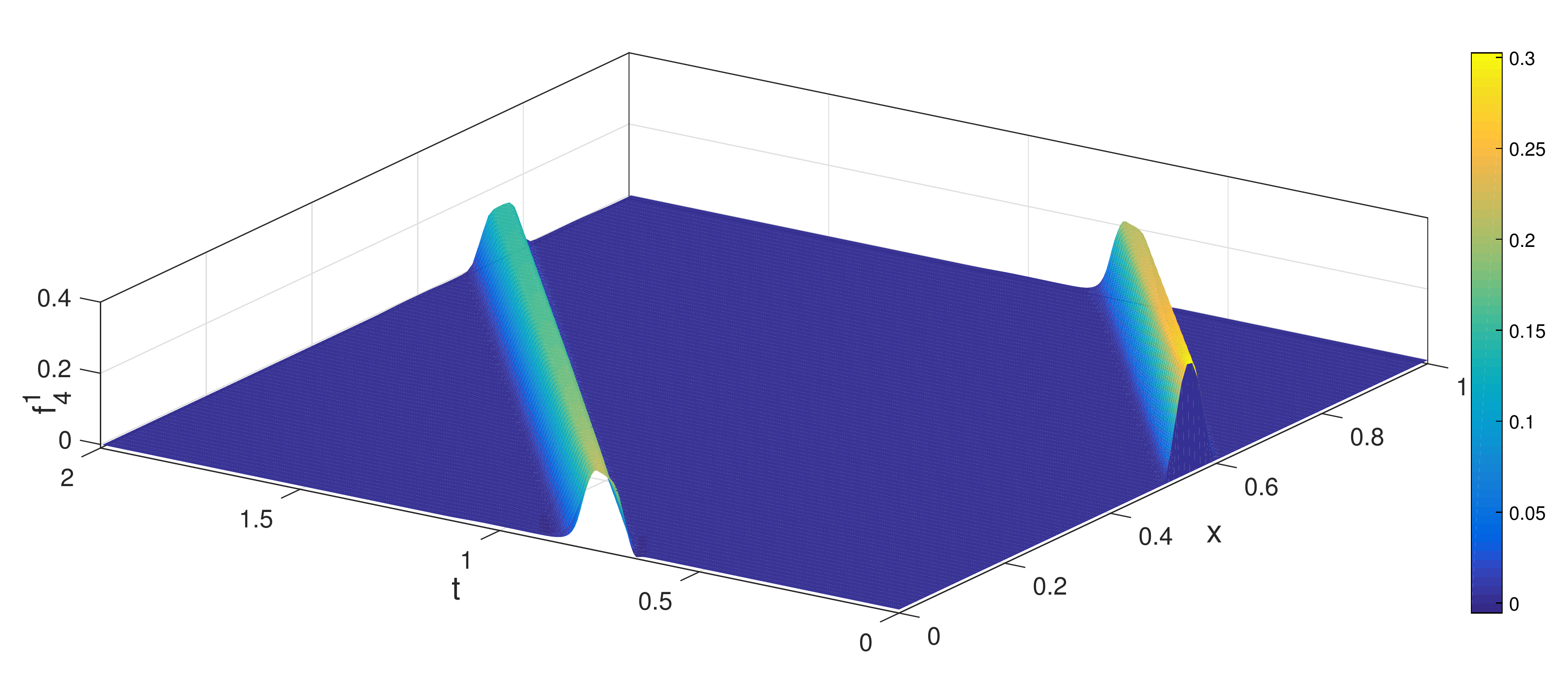}}
	~\subfigure[]{\includegraphics[height=1.5in ,width=2in]{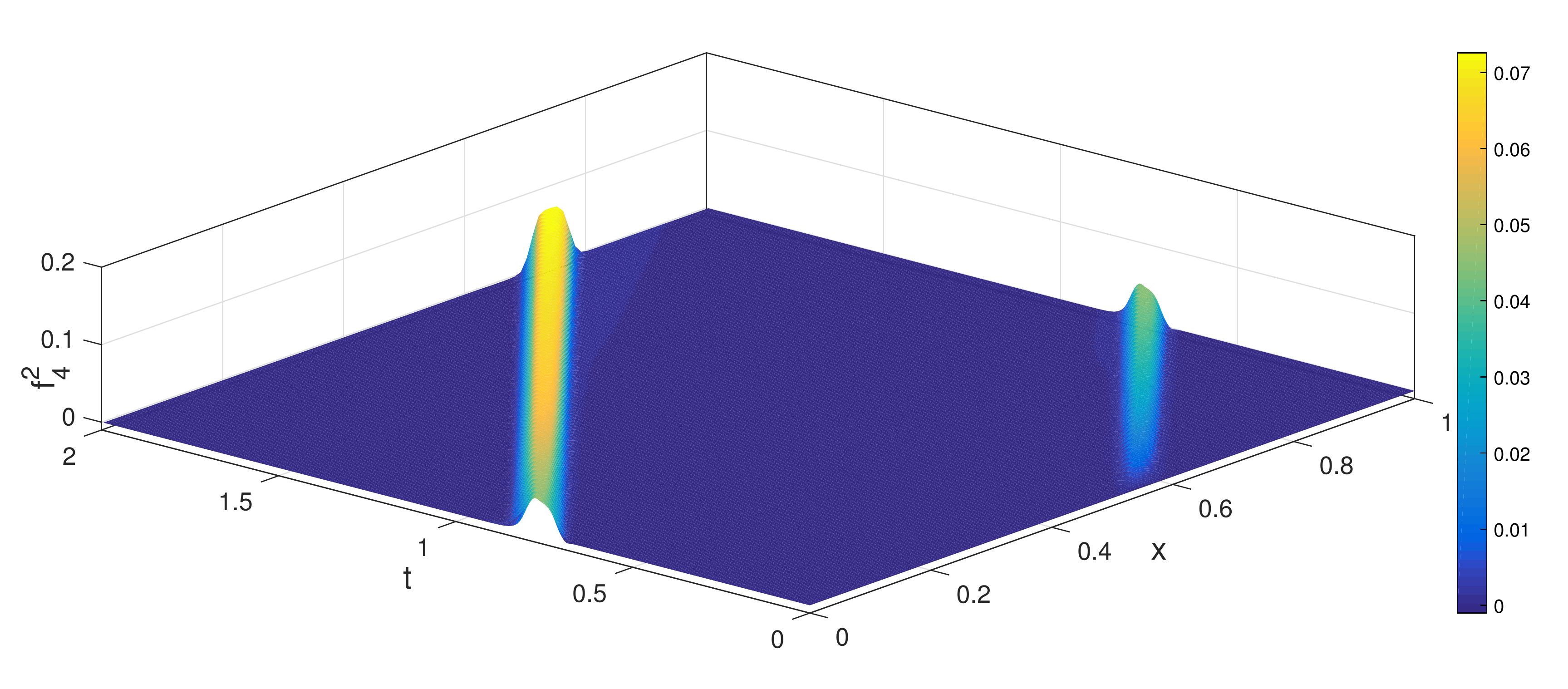}}
	\\{Evolution of $f_{4}^1$ (left) and $f_4^2$ (right).}\\
	\subfigure[]{\includegraphics[height=1.5in ,width=2in]{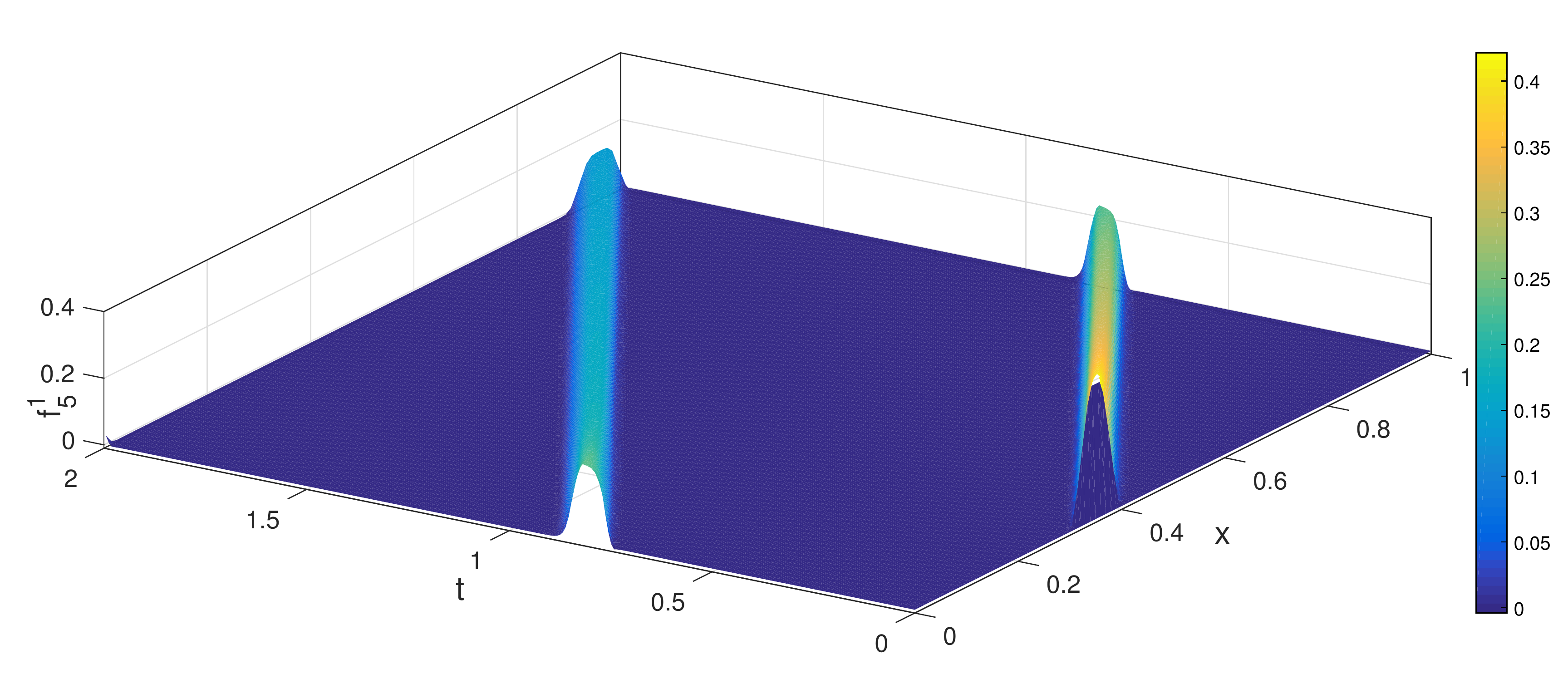}}
	~\subfigure[]{\includegraphics[height=1.5in ,width=2in]{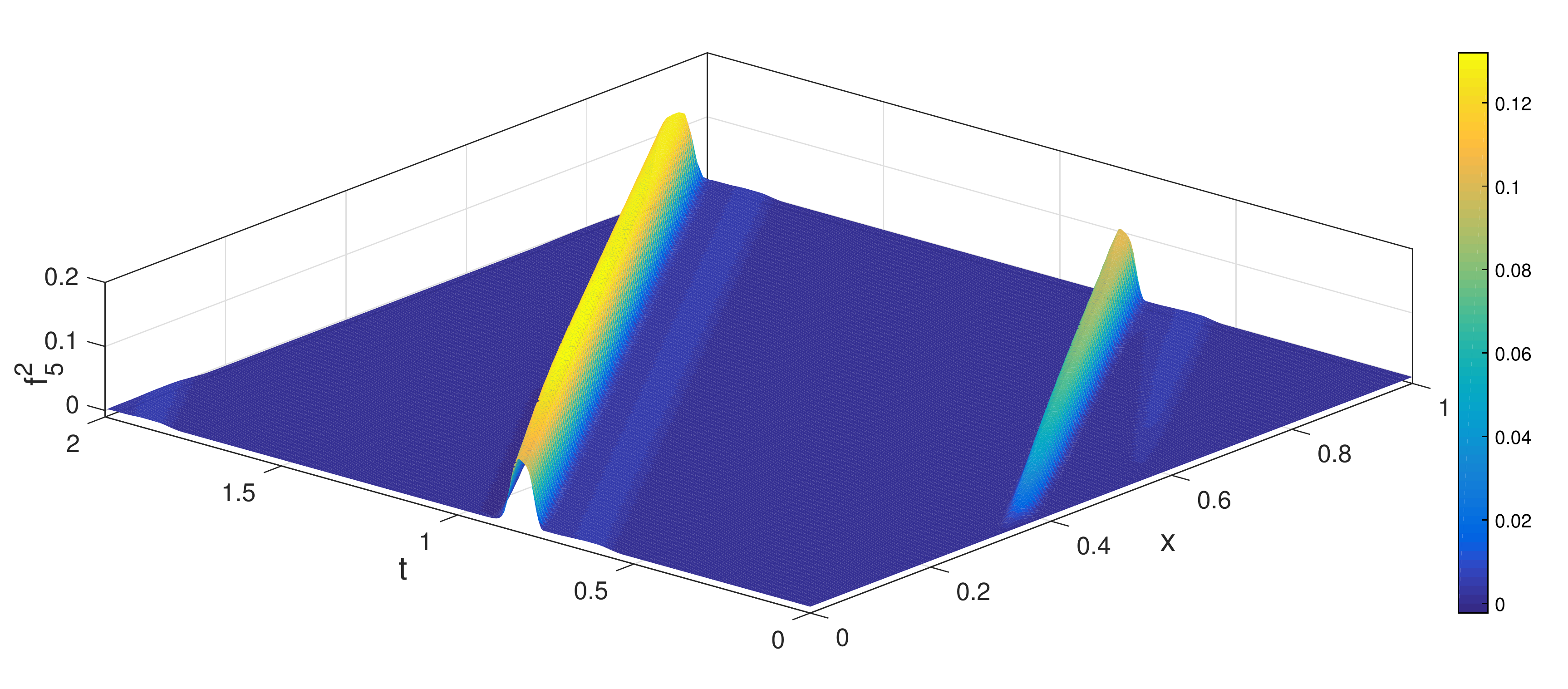}}
	\\{Evolution of $f_{5}^1$ (left) and $f_5^2$ (right).}\\
	\subfigure[]{\includegraphics[height=1.5in ,width=2in]{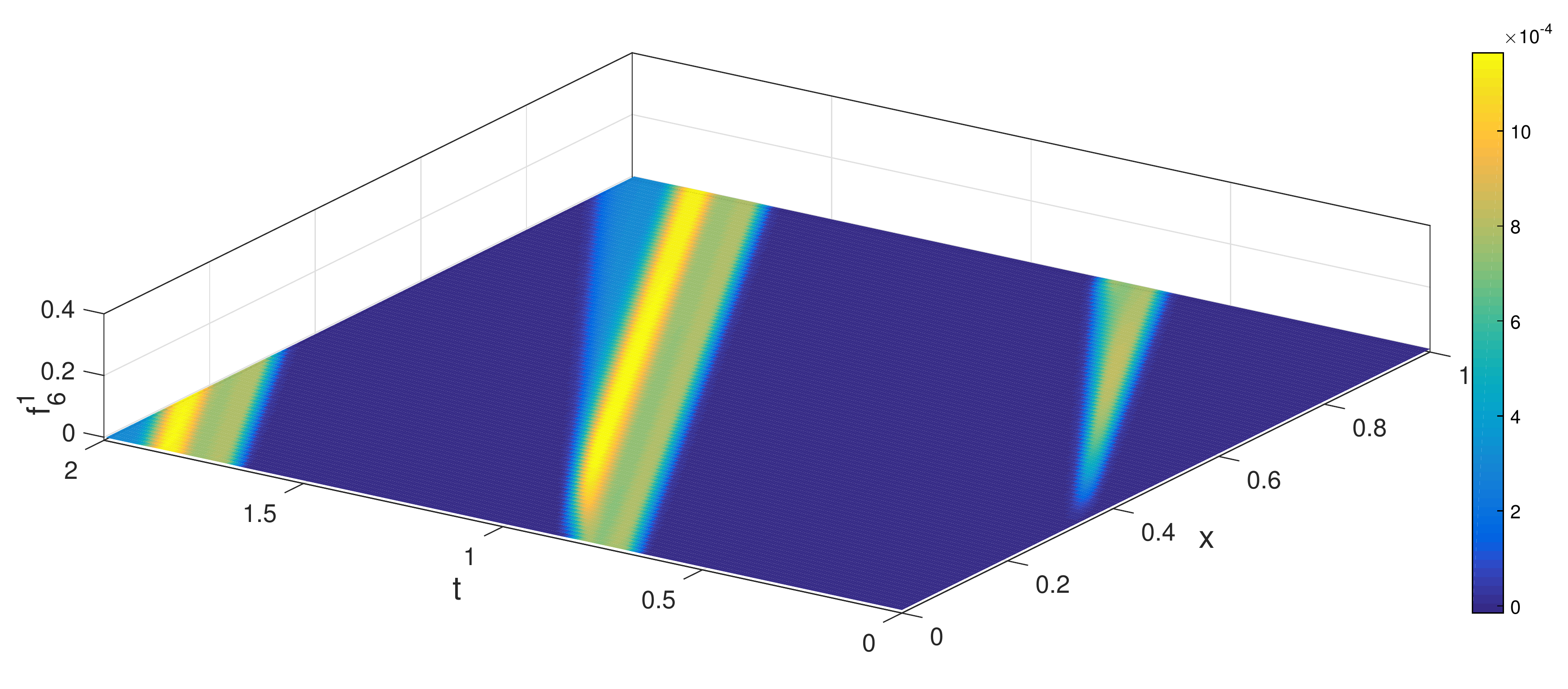}}
	~\subfigure[]{\includegraphics[height=1.5in ,width=2in]{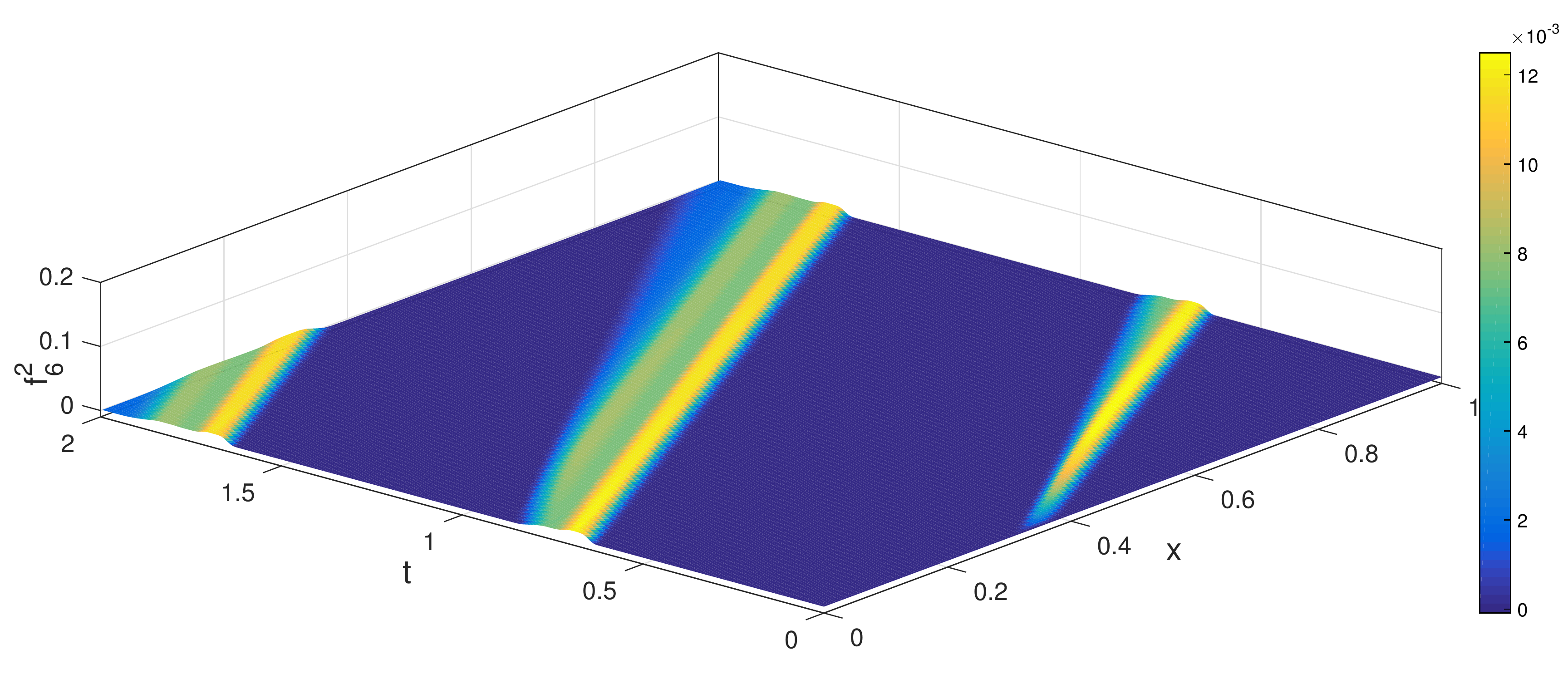}}
	\\{Evolution of $f_{6}^1$ (left) and $f_6^2$ (right).}\\
	\caption{Evolution of the micro-system in the in slowest lane (left) and in the fastest lane (right) under bad road conditions $\alpha=0.3$.}
	\label{1.8}
\end{figure}

Figure \ref{1.7} shows the evolution of groups of vehicles in the case of bad road conditions $\alpha=0.3$. We notice that:

\noindent \textbullet $\;$ The group of vehicles with velocity $v_{n-1}$ has approached the other in the slowest lane, as shown in Figure \ref{1.7}-(c), and then passes the other group of vehicles with the velocity $v_{n-2}$, see Figure \ref{1.7}-(e);\\
\noindent \textbullet $\;$ In the fastest lane, two groups of vehicles and a small group with velocity $v_n$ have appeared (see Figure \ref{1.8}-(d)), almost merging at a certain point. With the time being the group of vehicles with velocity $v_{n-1}$ passes the other group with velocity $v_{n-2}$ and the small group of vehicles takes the lead (see Figure \ref{1.7}-(f)). \\

The fist column in Figure \ref{1.8} shows the evolution time $t$ of the micro-systems $f_i^1$ for $i=3,\dots,6$ under the bad road conditions ($\alpha=0.3$). We notice that a negligible number of vehicles that would either decease (see Figure \ref{1.8}-(c) and \ref{1.8}-(e)) or increase (see Figure \ref{1.8}-(a) and \ref{1.8}-(g)). On the other hand, 
second column in Figure \ref{1.8} shows the evolution of the groups of vehicles in the fastest lane. We notice that the number of vehicles in the two groups $f_4^2$ and $f_5^2$ is increasing with the time being (see Figures \ref{1.7}-(d) and \ref{1.7}-(f)), and that a small group of vehicles with velocity $v_n$ have appeared (see Figure \ref{1.8}-(b) and \ref{1.8}-(h)).

\section{Conclusions and perspectives}
In this paper we have developed a new discrete kinetic model for multilane traffic flow on the basis of the kinetic theory approach. After a detailed description of a multilane road and the discrete method adopted, we have presented our proposed discrete model. Next, we have modeled each term appearing in the proposed model. We have proposed an improved table of games by taking into account further features of the complexity, namely the nonlinear additive interactions, perceived density rather than real one, and the road conditions as a function depending on the space variable which can used to depict the road maintenance areas. Moreover, we have modeled the external actions such as tollgates as well traffic signs. The well-posedness of the related Cauchy problem for the spatially homogeneous case has been proved by using Banach fixed-point theory. In order to validate the proposed improved model, in both the spatially homogeneous and inhomogeneous cases, and in fact to show its ability to reproduce certain realistic phenomena, we have reproduced the numerical simulations for Kerner's fundamental diagrams, and the asymptotic in time property along three lanes. Finally, we have obtained numerical simulations of the emerging of clusters with closed speed along two lanes by using the higher resolution finite volume method. We mention that all numerical simulations are provided with particular focus on the road environment conditions modeled by a parameter, which has been shown to have a great influence on the dynamic of vehicles.

Looking ahead, the obtained numerical simulations in the case of the spatially inhomogeneous problem are just one of the several emerging behaviors that can be depicted by our proposed modeling. One can also reproduce many realistic phenomena such as stop and go waves, and the external action by tollgates (see \cite{[Do1]} for one lane). In addition, it is interesting to develop an approach of  asymptotic limits towards macroscopic models known in the literature, for example, the recent models \cite{[GR20],[HMV18]}. One can draw inspiration from the ideas developed in \cite{[BD11]}.








\end{document}